\documentclass[10pt]{article}
\usepackage{graphicx}
\usepackage[breaklinks=true,colorlinks=true,linkcolor=blue,urlcolor=blue,citecolor=blue]{hyperref}
\usepackage{subfigure}
\usepackage{amsmath}
\usepackage{amsthm}
\usepackage{amssymb}
\usepackage{amsfonts}
\usepackage{tikz}

\begin{document}

\title{Cauchy difference priors for edge-preserving Bayesian inversion with an application to X-ray tomography}

\author{Markku Markkanen$^1$, 
Lassi Roininen$^2$, 
Janne M J Huttunen$^3$ \\ and 
Sari Lasanen$^4$\\ \\
$^1$ Eigenor Corporation \\
$^2$ University of Warwick, Department of Statistics \\
$^3$Nokia Technologies, and,  \\
University of Eastern Finland, 
Department of Applied Physics \\
$^4$ University of Oulu, Department of Mathematical Sciences
}

%
%
%
%

\maketitle

\begin{abstract}
We study Cauchy-distributed difference priors for edge-preserving Bayesian statistical inverse problems. 
On the contrary to the well-known total variation priors, one-dimensional Cauchy priors are non-Gaussian priors also in the discretization limit.
Cauchy priors have  independent and identically distributed increments.
One-dimensional Cauchy and Gaussian random walks are special cases of L\'evy $\alpha$-stable random walks with $\alpha=1$ and $\alpha=2$, respectively.
Both random walks can be written in closed-form, and as priors, they provide smoothing and edge-preserving properties.
We briefly  discuss also continuous and discrete L\'evy $\alpha$-stable random walks, and generalize the methodology to two-dimensional priors.

We apply the developed algorithm to one-dimensional deconvolution and two-dimensional X-ray tomography problems.
We compute conditional mean estimates with single-component Metropolis-Hastings and maximum a posteriori estimates with Gauss-Newton-type optimization method.
We compare the proposed tomography reconstruction method to filtered back-projection estimate and conditional mean estimates with Gaussian and total variation priors.
\end{abstract}

%
\vspace{2pc}

%
%
%

\section{Introduction}

In Bayesian statistical inverse problems, the objective is to estimate posterior distribution of an unknown object $\mathcal X$, given noise-perturbed indirect measurements. 
In this paper, we consider only linear inverse problems.
The measurements are formally given with an operator equation 
\begin{equation} \label{eqn:observationmodel}
m = \mathcal{AX} + e,
\end{equation}
where $m\in \mathbb{R}^n$ is a known measurement vector, and $\mathcal{A}$ is a known linear operator from some function space to a finite-dimensional vector space. 
We assume a Gaussian-distributed measurement noise $e\sim \mathcal{N}(0,\Sigma)$, where $\Sigma\in\mathbb{R}^{n\times n}$ is a known covariance matrix. We furthermore assume that $e$ is  independent of $\mathcal{X}$.

For computational purposes, we discretize the continuous observation model in Equation (\ref{eqn:observationmodel}) and denote the discretized observation model as
\begin{equation} \label{eqn:discreteobservationmodel}
m = AX + e,
\end{equation}
where  $A\in\mathbb{R}^{n\times k}$, and $X\in\mathbb{R}^k$.
A posteriori probability distribution is the solution of a finite-dimensional Bayesian statistical inverse problem.
Via the Bayes' formula, we give the posterior as an unnormalized probability density
\begin{equation} \label{eqn:bayes}
D(X|m) = \frac{D(X)D(m|X)}{D(m)} \propto	D(X)D(m|X).
\end{equation}
$D(X)$ is the a priori probability density, which reflects our information of the unknown before any actual measurement is done.
Prior is, in practice, the only tunable parameter in the estimation algorithm.
$D(m|X)$ is the likelihood density, which we construct from the discretized observation model (\ref{eqn:discreteobservationmodel}). 
$D(m)$ is simply a normalization constant, and hence we may drop it and use the unnormalized density.

In our previous papers (Roininen et al. 2011, 2013 and 2014 \cite{Roininen2011,Roininen2013,Roininen2014}), our main result has been that the prior should be modeled as a continuous-parameter random field, and for practical computations we discretize the prior in some computationally efficient way.
In the previous papers, we have considered Gaussian Markov random fields within the framework of Bayesian statistical inverse problems through, and we modeled the prior covariance via the sparse Cholesky decomposition of the inverse covariance matrix.
In \cite{Lindgren2011},  Lindgren et al.\ 2011 considered spatial interpolation  in a similar setting.
Such constructions promote discretization-invariant inversion, which means, in practice, that the posterior distributions are essentially the same, and hence also point estimators are  same on dense enough meshes.
For the considered Gaussian priors, it is enough to show that the discrete priors converge to continuous priors in the discretization limit, and more specifically, we only need to show that the discrete covariance matrix converges to a continuous covariance in the discretization limit. 
This would guarantee the convergence of the posterior distributions, and hence discretization-invariance, as shown by Lasanen 2012 \cite{Lasanen2012a,Lasanen2012b}.

Gaussian priors are known to smooth the edges in the reconstructions. 
Therefore, common alternative is to use total variation (TV) priors.
However, TV prior is not discretization-invariant as shown by Lassas and Siltanen 2004 \cite{Lassas2004}. This means that when we make the computational mesh denser and denser, the posterior distributions and estimates are not invariant with respect to the discretization of the unknown. 
For low-dimensional TV priors, we may get edge-preserving inversion, however,  high-dimensional TV priors start to loose their edge-preserving property when  the  random walk presenting the prior  tends to its  continuous limit.
From  Donsker's theorem, Lindeberg-L\'evy theorem and  Doeblin's theorem,  we see that  the problem originates from the  distributions of  the random walk. 
Namely,  the jumps have 
finite  moments, which always leads to normally distributed limits. When one allows other distributions for jumps, it may  become possible to maintain the edge-preserving properties for high-dimensional priors.  

Lassas, Saksman and Siltanen 2009 \cite{Lassas2009} proposed to use non-Gaussian Besov space priors which are constructed on wavelet basis.
These can be shown to promote discretization-invariant inversion, and the method has been e.g.\ applied to X-ray tomography by V\"ansk\"a et al.\ 2009 \cite{Vanska2009}.
In this paper, we shall tackle the issue of edge-preserving inversion  by constructing Cauchy difference priors starting from continuous Cauchy random walks. 
We apply the methodology to deconvolution and X-ray tomography problems.
Our discussion lies mostly on the numerics and we leave mathematically rigorous discussion on discretization-invariance to the future papers.

The rest of this paper is organized as follows: 
In Section \ref{sec:onedimensional}, we consider one-dimensional Cauchy priors. We start by considering L\'evy $\alpha$-stable random walks and consider why Cauchy random walks promote edge-preserving inversion. We then review single-component Metropolis-Hastings algorithm, and use the developed methodology to a one-dimensional deconvolution problem. 
In Section \ref{sec:twodimensional}, we consider two-dimensional Cauchy  priors and apply the methodology to a two-dimensional X-ray tomography. We compare the developed algorithm against conditional mean estimates with Gaussian and TV priors as well as filtered back-projection reconstructions. 
In Section \ref{sec:conclusion}, we conclude the study and make some notes on future research.

\section{One-dimensional Cauchy priors in Bayesian inversion}
\label{sec:onedimensional}
In this section, we first consider L\'evy $\alpha$-stable random walks, and obtain the Cauchy and Gaussian random walks as special cases.
We study  edge-promoting  inversion idea with Cauchy priors  via unimodality versus bimodality of the distributions derived from the Cauchy distributions in a specific way. 
We carry out a comparison with other common priors.
Finally, we apply the methodology to one-dimensional deconvolution problem by obtaining conditional mean estimates with single-component Metropolis-Hastings.

\subsection{Defining random walks}

Let $\{\mathcal X(t) ,t\in \mathbb{I}\subset \mathbb{R}^+\}$ be a stochastic process. 
We call it a continuous L\'evy $\alpha$-stable process starting from zero, if $\mathcal X(0)=0$, $\mathcal X$ has independent increments and  
\begin{equation}
\mathcal X(t)-\mathcal X(s) \sim S_\alpha\left(\left(t-s\right)^{1/\alpha}, \beta,0\right)
\end{equation}
for any $0\leq s < t < \infty$ and for some $0<\alpha\leq 2, -1\leq \beta \leq 1$ 
(for stable distributions, see e.g.\ \cite{MR1280932}). 

We can simply construct stable processes by using independently scattered measures.
A random measure $M:\Omega\rightarrow 
L^0(P)$  is called independently scattered if $M(A_1),\dots,M(A_n)$ are 
independent random variables whenever the sets $A_i$, $i=1,\dots,n$, $n\in\mathbb N$,
are pairwise disjoint.  We require that  
\begin{equation}
M(A)\sim S_\alpha \left( \vert A \vert^\frac{1}{\alpha}, \beta,0\right)
\end{equation} 
with constant skewness $\beta$. Then $\alpha$-stable L\'evy motion starting from zero is
\begin{equation}
\mathcal X (t) = M\left(\chi([0,t])\right),
\end{equation}
where $\chi_A$ is the characteristic function of the set $A$.

For the continuous limit of the L\'evy $\alpha$-stable random walk, we apply independently scattered measures. 
Let us denote the discrete random walk at $t=jh$ by $X_{j}$, where $j\in \mathbb{Z}^+$ and $h>0$ is the discretization step. 
We obtain a random walk approximation (See Chapter 3.3. in \cite{MR1280932})
\begin{equation} \label{eqn:Levywalk}
X_{j}-X_{j-1} \sim S_\alpha \left( h^\frac{1}{\alpha}, \beta,0\right).
\end{equation}
It is easy to see that 
such a random walk approximation converges to the continuous L\'evy $\alpha$-stable  random walk as $h\rightarrow 0$ in distribution in the Skorokhod space of  functions that are right-continuous and have left limits.

%
%
%
%
%
%

Our special interest is the Cauchy random walk.
Given Equation (\ref{eqn:Levywalk}) and choosing $\alpha=1$, $\beta=0$, we can write the Cauchy random walk increments as a probability  density
\begin{equation}
D (X) = C \prod_{j=1}^J \left( \frac{\lambda_j h}{(\lambda_j h)^2 + (X_j-X_{j-1})^2} \right),
\end{equation}
where $\lambda_j$ is the regularizing term. 
For the sake of simplicity, let us assume $\lambda_j = \lambda$ to be a constant.
$C$ is a normalization constant.
We can apply this probability density directly as an a priori probability density in Bayes' formula in Equation (\ref{eqn:bayes}).

\begin{figure}[htp]
\begin{center}
  \subfigure[Cauchy random walk]{\includegraphics[width=0.49\textwidth]{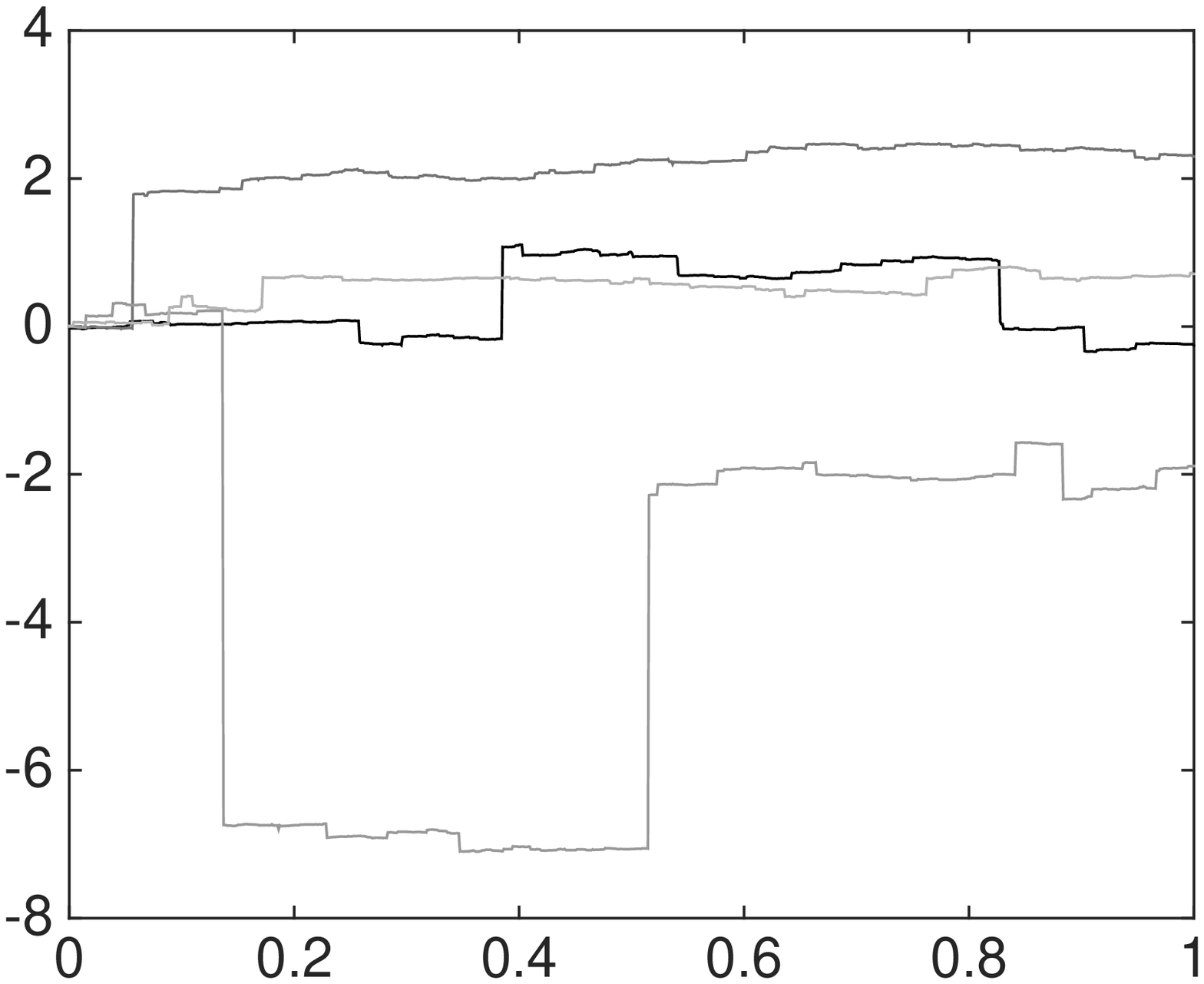}}
  \subfigure[Gaussian random walk]{\includegraphics[width=0.49\textwidth]{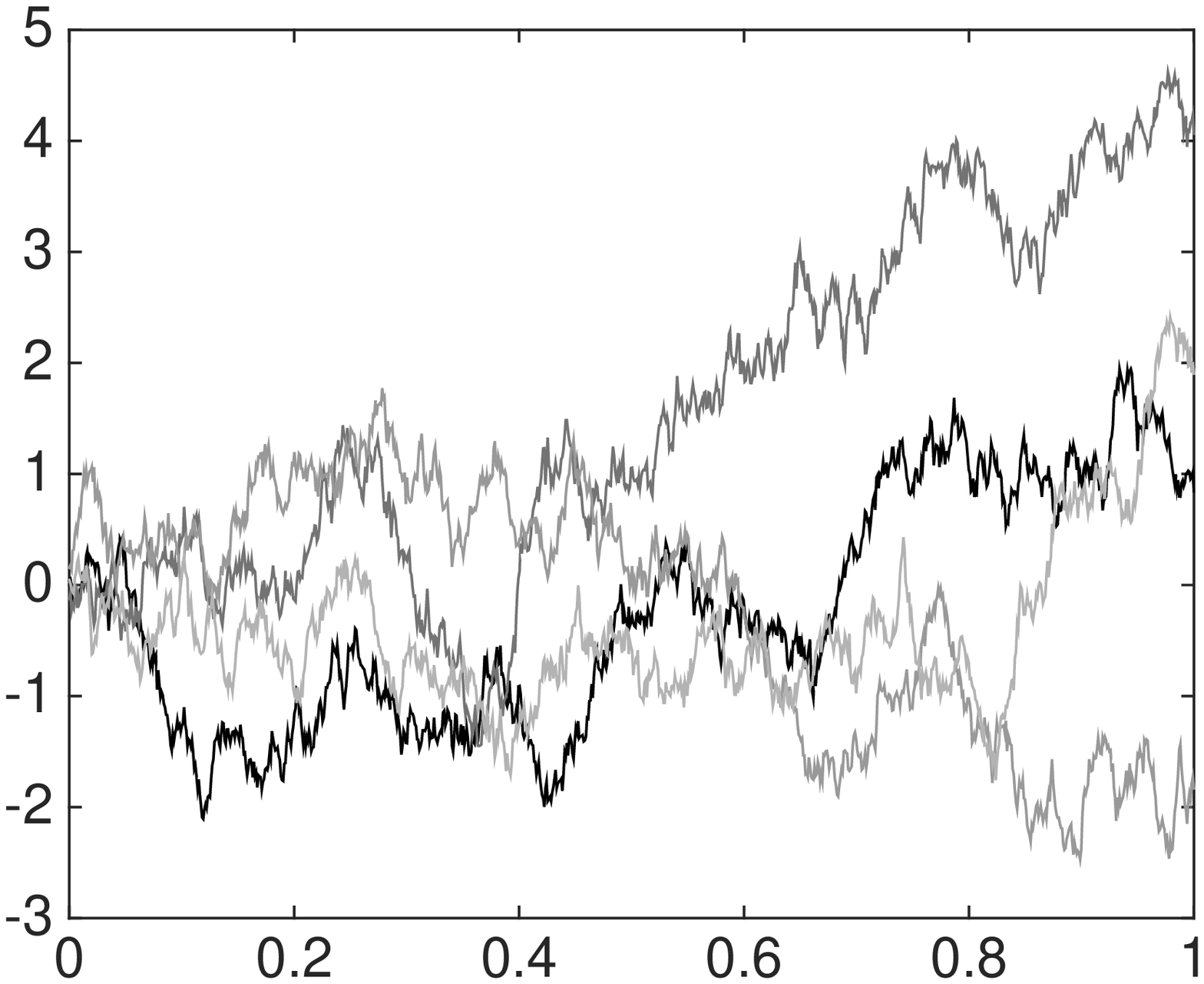}}

  \caption{Realizations of Cauchy and Gaussian random walks.} \label{fig:random_walks}
  \end{center}
\end{figure}

In Figure \ref{fig:random_walks}, we have plotted  realizations of Cauchy and Gaussian random walks. The latter is the L\'evy $\alpha$-stable random walk with $\alpha = 2$.
The Cauchy random walk realizations have either small perturbations or big jumps. 
%
Hence, intuitively it feels plausible that using Cauchy distributions for differences lead to edge-preserving inversion.
%
%
Gaussian random walk realizations do not prefer jumps, but continuous paths, hence the tendency for smoothing in Bayesian inversion.

\subsection{Edge-preserving inversion}
%

To our best knowledge, for difference priors there has not been any simple criteria which would indicate, whether the prior favors edge-preserving solutions instead of smooth ones.
%
The intuitive idea is to construct a prior, which promotes discontinuities, i.e.\ jumps. 
Hence, via a negation, we may say that in edge-preserving inversion, the priors do not promote smooth properties.

\begin{figure}[htp]
\begin{center}
  \subfigure[Unimodal Cauchy $a<1$]{\includegraphics[width=0.32\textwidth]{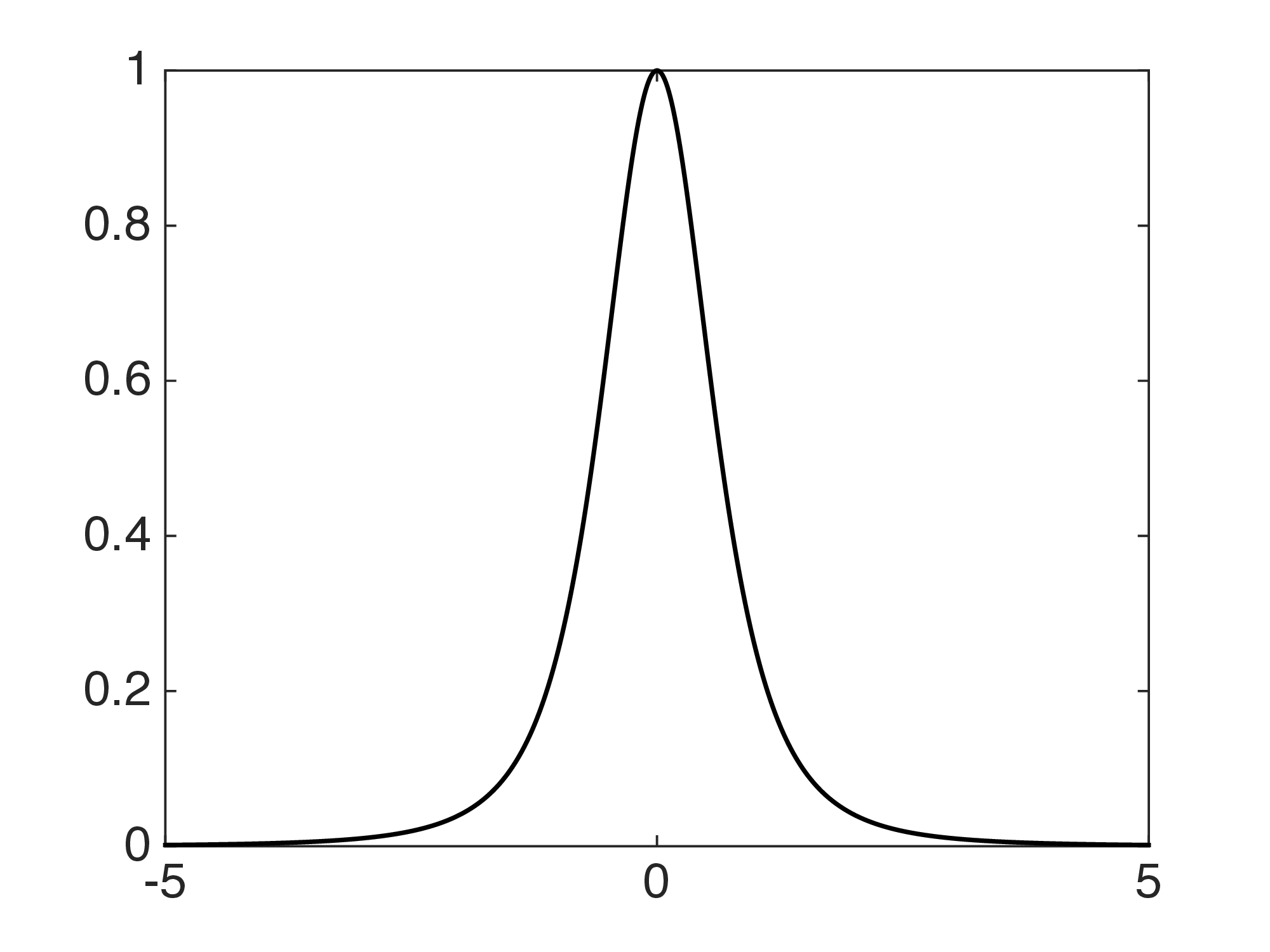}}
  \subfigure[Flat Cauchy with $\vert a\vert=1$]{\includegraphics[width=0.32\textwidth]{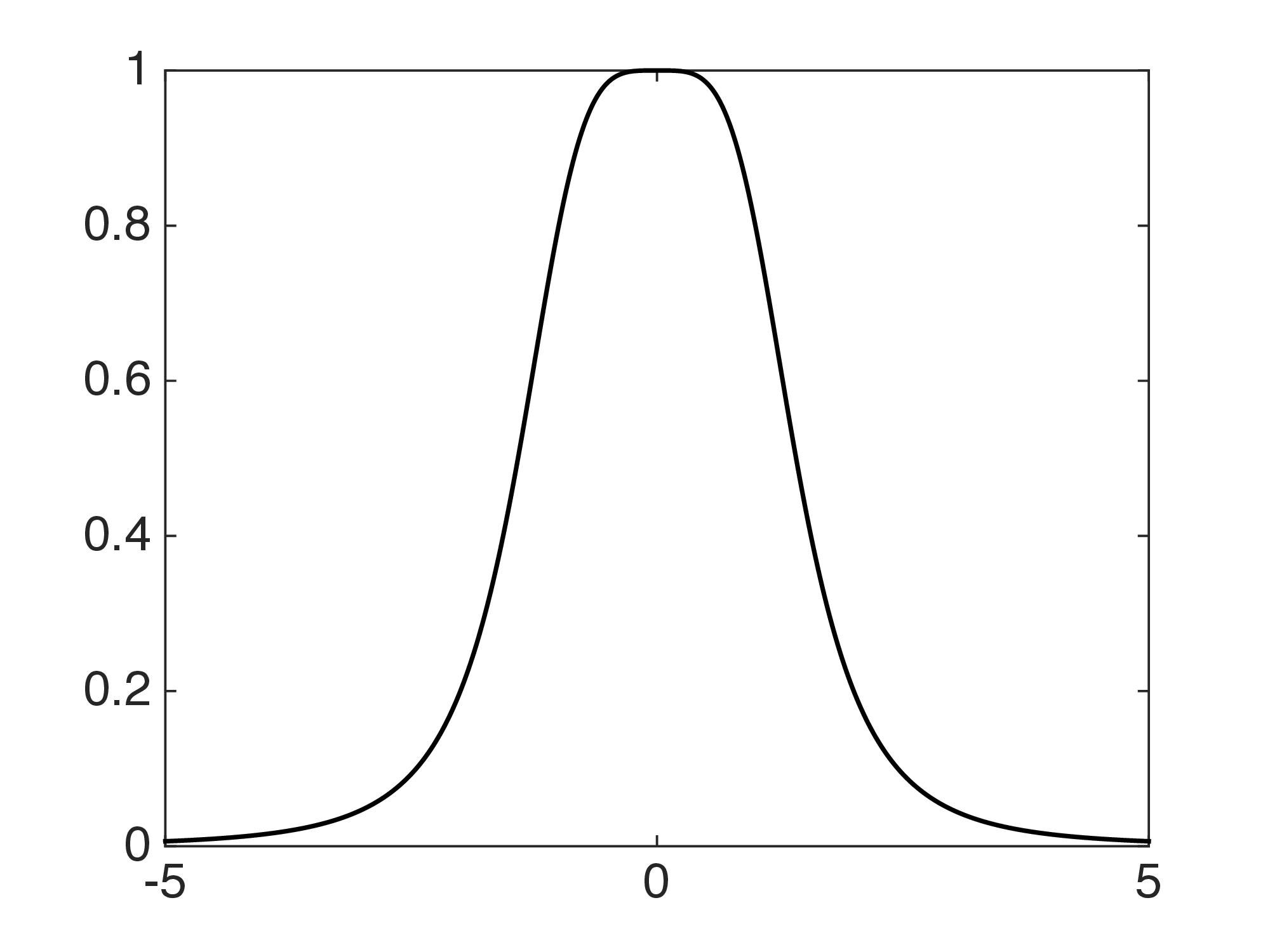}}
  \subfigure[Bimodal Cauchy $a >1$]{\includegraphics[width=0.32\textwidth]{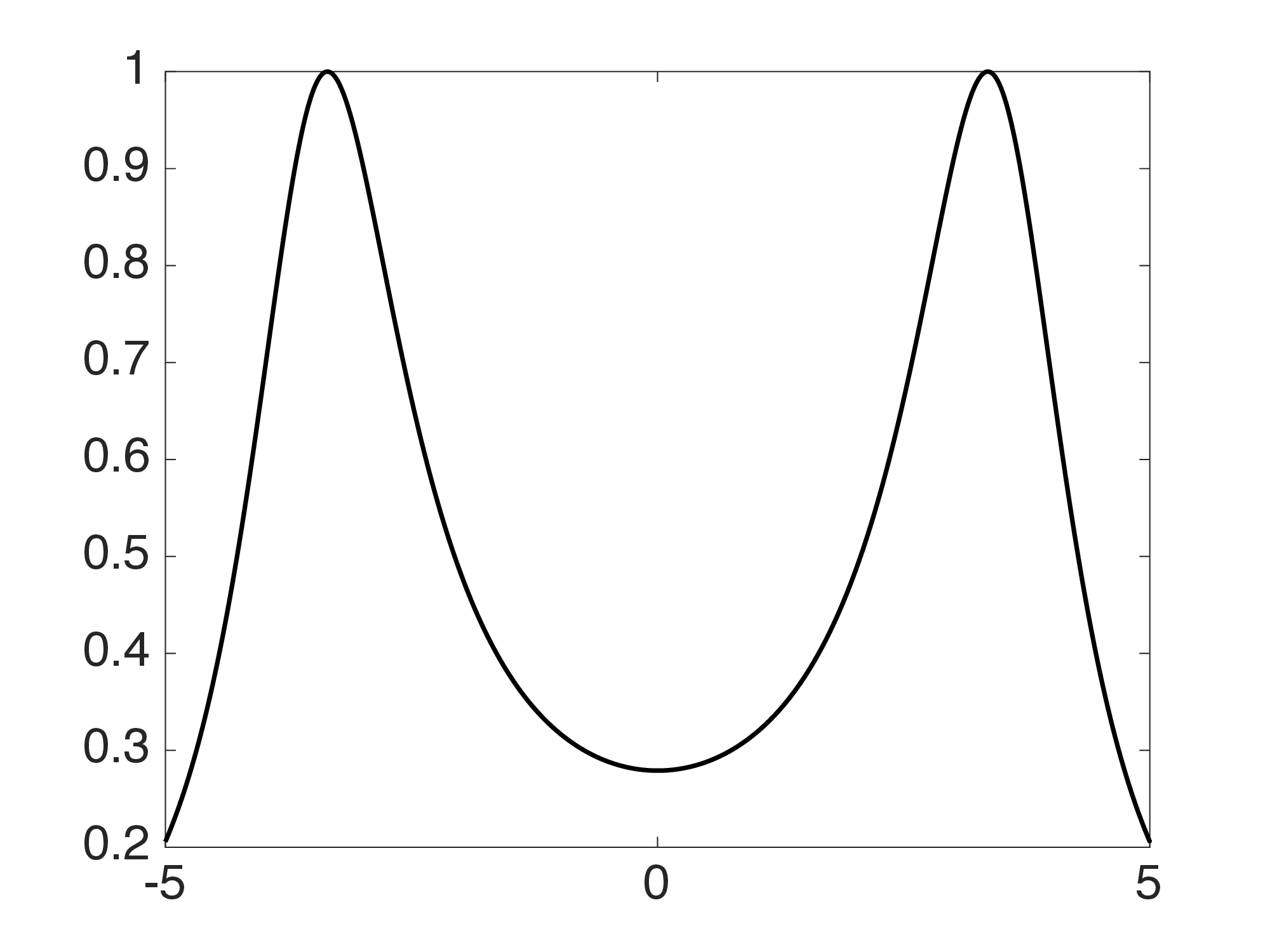}}

  \subfigure[Gaussian]{\includegraphics[width=0.32\textwidth]{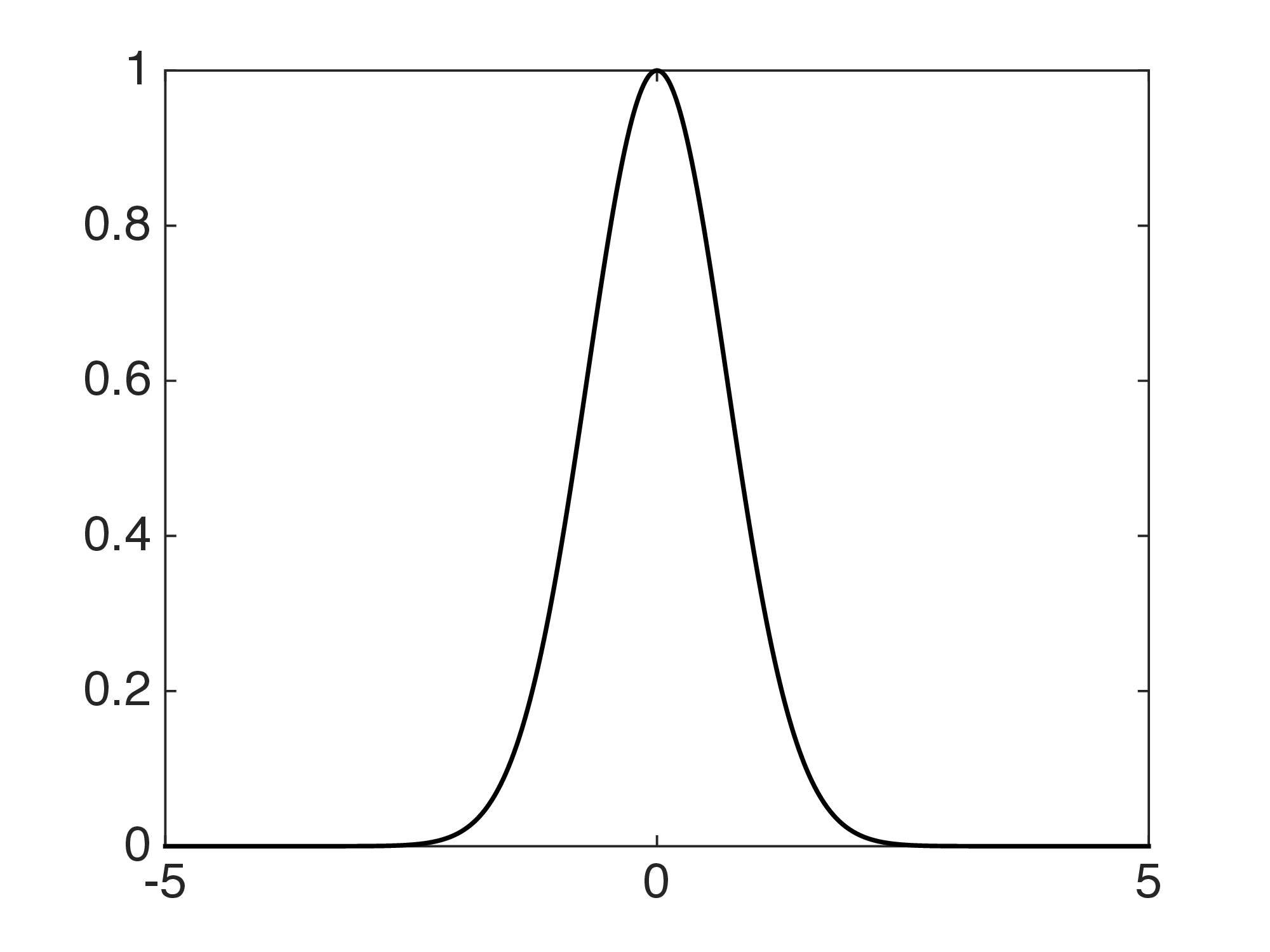}}
  \subfigure[Total variation]{\includegraphics[width=0.32\textwidth]{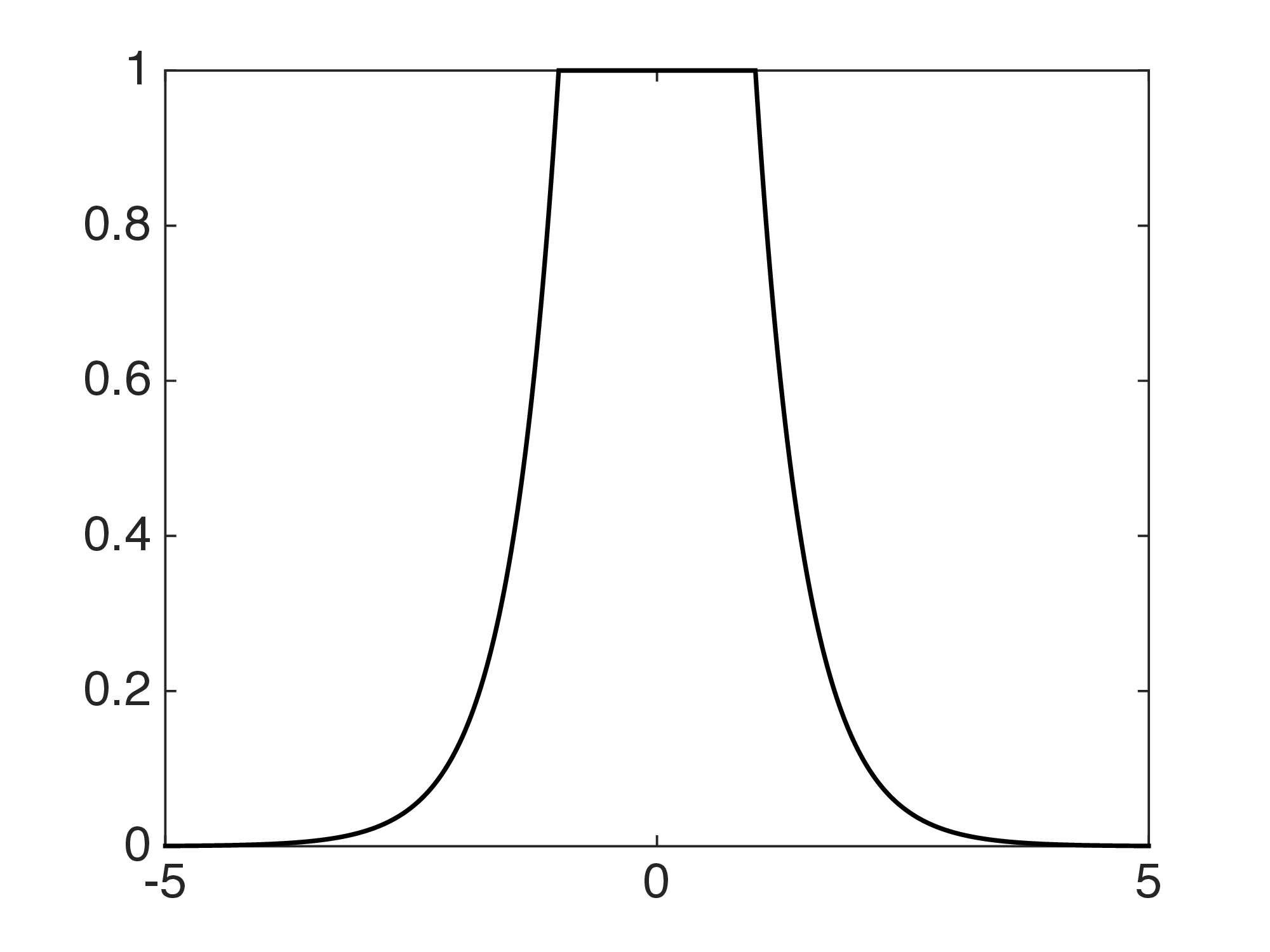}}  

  \caption{Upper panel: Cauchy probability density function for $X_2$,
    given fixed $X_1=-a$ and $X_3=a$. Bottom panel: The same case for
    Gaussian difference prior and total variation
    prior.
  } \label{fig:probability_densities}
  \end{center}
\end{figure}

Here, we will  consider a one-dimensional case, as the situation becomes more complicated in higher dimensions, and look at three consecutive coordinates $X_{j-1},X_j,X_{j+1}$ of the random vector $X$. 
We could then say that $X$ prefers smoothness at point $j$, if $X_j$ is (most probably) in the middle of  $X_{j-1}$ and $X_{j+1}$, and that there is edge at point $j$ if  $X_j$ is (most probably) either at $X_{j-1}$ or at $X_{j+1}$. 
The above mentioned idea means that we are interested in the conditional prior distribution of $X_j$ on the condition that $X_{j-1}$ and $X_{j+1}$ are known. 
Without loss of generality we can then assume that $X_{j-1} = -a$ and $X_{j+1} = a$. 
%
%
%
%
%
The Cauchy difference prior for $X_j$ can then be written as a probability density
\begin{equation}\label{eq:1}
D(X_j) \propto \frac{1}{1+(X_j-a)^2}\frac{1}{1+(X_j+a)^2}.
\end{equation}
This is simply a product of two Cauchy probability density functions, one from each neighbor of $X_j$.
In order to see the properties of these functions, consider the second derivative of the probability density function with respect to $X_j$ at zero:
\begin{equation}
\begin{split}
D''(0) < 0, \textrm{when},~|a| < 1~:~\textrm{maximum at}~0, ~\textrm{unimodal} \cr
D''(0) = 0,  \textrm{when},~|a| = 1~:~\textrm{as flat as possible at}~0 \cr
D''(0) > 0,  \textrm{when},~|a| > 1~:~\textrm{minimum at}~0,~\textrm{bimodal}.%
\end{split}
\end{equation}
In Figure \ref{fig:probability_densities}, we have plotted these probability density functions. 
With $|a| < 1$ we have a unimodal function, hence regularization promotes smooth solution around $X_j=0$.
With $|a|>1$ we have a bimodal function, hence regularization promotes solutions around $X_j = \pm a$ with smooth variations around these two values. 
This means that this kind of prior promotes jumps or small local oscillations, and this is how we understand edge-preserving prior in this paper.

We can make similar plots for Gaussian difference prior and  total variation prior. 
The Gaussian difference prior is simply a Gaussian-shaped function with maximum at $X_j=0$. 
Hence, providing smoothing properties.
The total variation priors is constant for $[-a,a]$ and decays exponentially outside this region.
Hence, total variation neither punishes discontinuities or promotes edge-preserving inversion as it promotes all values equally between these two points.
This is also the case for the Cauchy prior with $|a|=1$, which is, in a sense, similar to the total variation prior.

%
%

\subsection{Single-component Metropolis-Hastings}

In order to get estimators for Bayesian statistical inverse problems with Cauchy difference priors, we need to compute either the maximum a posteriori (MAP) estimate or conditional mean (CM) estimates.
For the MAP estimation, we use Gauss-Newton type optimization methods. 
%
However, our interest lies mostly in the  CM estimate, which is the mean of the posterior density (see Equation (\ref{eqn:bayes})).
In higher dimensional problems, the explicit computation of the CM
estimates, the integration over the posterior density, is rarely
possible and numerical integration of posterior density typically
leads to an infeasible approach.
Therefore, CM estimates are usually computed using Markov Chain Monte Carlo (MCMC)
sampling methods, such as Metropolis-Hastings (MH) algorithm or Gibbs
sampler (see e.g.\ \cite{Gilks,gelman03}).
The idea of the MCMC methods is to generate samples $\left\{X^{(s)}\right\}$
from the posterior distribution and approximate the CM estimate using
the ensemble mean.

In this paper, we test the proposed approach using one--dimensional deconvolution problem and
X-ray tomography problem. 
Both of these problems have a specific property that, when evaluating the
posterior density, an update of a single component in the parameter
vector is computationally very cheap.
The commonly known algorithms that take advantage of the property are
single-component Gibbs sampler and  single-component
Metropolis-Hastings (MH)
algorithm; see for example \cite{Gilks}. 
%
Disadvantage of the Gibbs sampler is that  the explicit form of the
conditional distributions is known, but the inverse cumulative is not known, and that would be needed for random number generation.

The idea of the single-component MH sampling is similar to the single
component Gibbs sample, except the component-wise sampling is carried
out using MH-type sampling instead of the direct sampling from the
conditional distributions.
We start with some initial value
\begin{equation}
  X^{(0)}=\left(X_1^{(0)},X_2^{(0)},\ldots,X_n^{(0)}\right)
\end{equation}
and set $s=0$.
First, we consider one--dimensional conditional
distribution of $X_1$ given data $m$ and all other parameters:
\begin{equation}
\label{eq:marginal_mu}
D\left(X_1|m,X_2^{(0)},X_3^{(0)},\ldots,X_n^{(0)}\right)
=C_1   D\left(X_1,X_2^{(0)},X_3^{(0)},\ldots,X_n^{(0)} \vert m\right).
\end{equation}
where $C_1$ is a normalization constant. 
We apply the MH algorithm to this distribution: We draw a candidate
sample $\tilde X_1$ from a proposal distribution $Q_1\left(\cdot|X_1^{(0)}\right)$,
which is accepted with a probability
\begin{equation}
  \alpha_1=\min\left(1,\frac{D\left(\tilde X_1|m,X_2^{(0)},\ldots,X_n^{(0)}\right)q_1\left(X_1^{(0)}|\tilde
      X_1\right)}{D\left(X_1^{(0)}|m,X_2^{(0)},\ldots,X_n^{(0)}\right)q_1\left(\tilde X_1|X_1^{(0)}\right)}\right).
\end{equation}
If accepted, we set $X_1^{(1)}=\tilde X_1$. Otherwise we set $X_1^{(1)}=X_1^{(0)}$.
We continue to the next parameter and consider the distribution
\begin{equation}
\label{eq:marginal_mu}
D\left(X_2|m,X_1^{(1)},X_3^{(0)},\ldots,X_n^{(0)}\right)
=C_2   D\left(X_1^{(1)},X_2,X_3^{(0)},\ldots,X_n^{(0)}|m\right).
\end{equation}
where $C_2$ is a normalization constant. 
We draw a candidate
sample $\tilde X_2$ from a proposal distribution
$q_2(\cdot|X_2^{(0)})$ and accept the sample with probability
\begin{equation}
  \alpha_2=\min \left(1,\frac{D\left(\tilde X_2|m,X_1^{(1)},X_3^{(0)},\ldots,X_n^{(0)}\right)q_2\left(X_2^{(0)}|\tilde
      X_2\right)}{D\left(X_2^{(0)}|m,X_1^{(1)},X_3^{(0)},\ldots,X_n^{(0)}\right)q_2\left(\tilde X_2|X_2^{(0)}\right)}\right).
\end{equation}
The procedure is continued until all parameters have been updated one at a time and we have
$X^{(s+1)}$. 
We increase $s\leftarrow s+1$ and repeat the sample generation $N_s$
times to obtain a set of samples $\left\{X^{(s)}\right\}_{s=1}^{N_s}$.
See e.g.\ \cite{Gilks} for more details.
The algorithm is called also  Metropolis-within-Gibbs \cite{Dunlop2016}.

In this paper the proposal distribution is chosen to be Gaussian:
$Q_i\left(\cdot|X_i^{(0)}\right)= \mathcal{N}(X_i^{(0)},\sigma^2$).
The variance $\sigma^2$ is chosen such that 25-50\% of samples are accepted.
It is common to ignore a number of samples at the beginning (burn-in
period) due to the fact that it takes a time to converge to the stationary
distribution $D(X|m)$. 
In this paper the burn-in period is chosen to be the half of the
sample size.
We leave optimization of the burn-in period to later studies.



\subsection{Numerical example: One-dimensional deconvolution}

%

\begin{figure}[htp]
\begin{center}
\includegraphics[width=\textwidth]{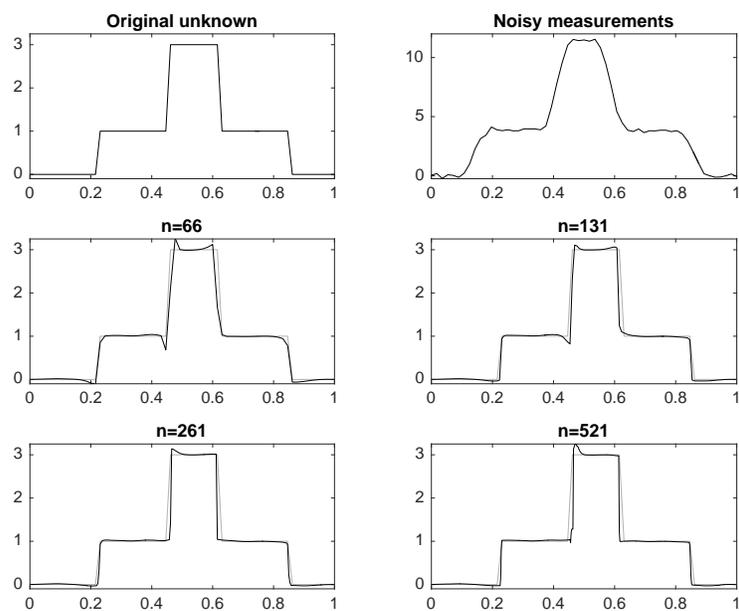}
  \caption{Top left panel: The original unknown. Top right panel: Noisy deconvolved measurement vector. Middle and bottom panels: Reconstructions with Cauchy priors on different grids with $n=66,~131,~261,~521$.
  } \label{fig:deconvolution}
  \end{center}
 \end{figure}

In order to demonstrate the edge-preserving and discretization-invariance properties, in Figure \ref{fig:deconvolution}, we have a synthetic deconvolution example with Cauchy  priors on different grids.
The unknown is assumed to be a piecewise constant function (upper left panel) and it is convolved with a constant cosine-shaped convolution kernel.
In the measurements, we have used additive white noise with 1\% noise level.

We have made the reconstructions with number of unknowns being $n=66,~131,~261,~521$. 
In the single-component Metropolis-Hastings, we use 1,000,000 long chain.
In order to save memory, we use only every 10th sample (as 10 consecutive samples are still nearly the same).
This is because we want the samples to be essentially independent.
Burn-in period is chosen to be 500,000. 
The MCMC chain could be optimized, but as in this paper our goal is to simply demonstrate the Cauchy prior, we postpone the MCMC optimization studies for subsequent papers.

The reconstructions on different grids clearly provide edge-preserving features. 
The shapes vary between the different reconstructions, but this is mostly due to the MCMC chain, and, characteristically they are all rather similar.
%
This kind of behavior of estimators cannot be obtained with total variation priors, as they are not discretization-invariant prior.
For the total variation prior examples, see Lassas and Siltanen 2004 \cite{Lassas2004}.

\section{Application to two-dimensional X-ray tomography} 
\label{sec:twodimensional}

Let us consider a two-dimensional lattice $(hj,h'j')$, where $(j,j')\in \mathbb{I}^2\subset \mathbb{Z}^2$ and discretization steps to each coordinate directions $h,h'>0$.
A Gaussian difference prior can be constructed via difference equations of type
\begin{equation} \label{eqn:gaussdiffs}
\begin{split}
 X_{j,j'} - X_{j-1,j'} \sim \mathcal{N}(0,\sigma^2h/h') \cr
X_{j,j'} - X_{j,j'-1} \sim \mathcal{N}(0,\sigma^2h'/h).
\end{split}
\end{equation}
The term $\sigma^2>0$ is the regularization parameter. 
The prior is then constructed through products of conditional normal distributions
\begin{equation}
D_{\mathrm{pr}}(X) \propto \exp\left(-  (X^TL_1^TC_1^{-1}L_1X + X^TL_2^TC_2^{-1}L_2X)\right),
\end{equation}
where $L_1$ and $L_2$ are difference matrices to two different coordinate directions, and, $C_1=\sigma^2h/h' I$ and $C_2=\sigma^2h'/h I$, where $I$ is an identity matrix.
These constructions are rather well-known, see for example \cite{Lasanen2005}.
We emphasize that Equation (\ref{eqn:gaussdiffs}) is not explicitely solvable.

Let us consider with a simple example why we scale the variances by $h/h'$ and $h'/h$.
It is enough to consider only one difference equation, as the other one is obtained through same arguments.
We consider the grid densification  in two steps: 1) We consider denser discretization along the 'regularization coordinate direction', and, 2) Denser discretization along the 'other coordinate direction'.
Let $X_{j,j'} - X_{j-1,j'} \sim \mathcal{N}(0,\sigma^2)$. 
If we densify the lattice along the regularization coordinate direction $j$, this would correspond to one-dimensional Gaussian random walk and hence we need to add the discretization step in the variance, hence $X_{j,j'} - X_{j-1,j'} \sim \mathcal{N}(0,\sigma^2 h)$.

\begin{figure}[htp]
\begin{center}
\begin{tikzpicture}
	\draw (-5, -1.5) -- (-5, 1.5);
	\draw (-3.5, -1.5) -- (-3.5, 1.5);
	\draw (-2, -1.5) -- (-2, 1.5);
	\draw (-5,-1.5) -- (-2,-1.5);
	\draw (-5,1.5) -- (-2,1.5);

	\draw (5, -1.5) -- (5, 1.5);
	\draw (3.5, -1.5) -- (3.5, 1.5);
	\draw (2, -1.5) -- (2, 1.5);
	\draw (5,1.5) -- (2,1.5);
	\draw (5,-1.5) -- (2,-1.5);
	\draw (5,0) -- (2,0);
	\draw (-4.25,0) node {$X_1$};
	\draw (-4.25+1.5,0) node {$X_2$};
	\draw (4.25,0.75) node {$X_{22}$};
	\draw (4.25-1.5,-0.75) node {$X_{11}$};	
	\draw (4.25-1.5,0.75) node {$X_{12}$};	
	\draw (4.25,-0.75) node {$X_{21}$};	
	
	\draw [->] (-1,0) to (1,0) ;
	\draw (0,0.5) node {Grid densification};
\end{tikzpicture}
\end{center}
\caption{Grid densification to the second coordinate direction.}
\label{fig:dens}
\end{figure}
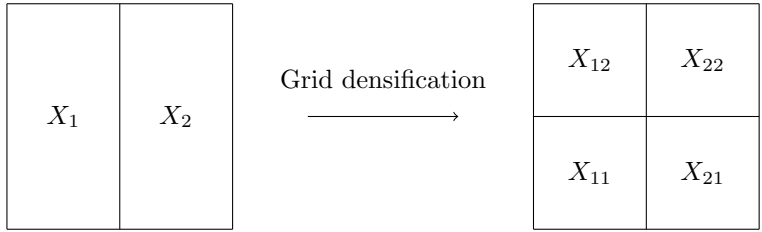

The second coordinate direction is a bit trickier.
Let us consider a simplified situation as depicted in Figure \ref{fig:dens}. 
Thus, we have
\begin{equation}
				X_1 = \frac{X_{11}+X_{12}}{2},~~\textrm{and},~~
				X_2 = \frac{X_{21}+X_{22}}{2}.
\end{equation}
Now the increments
\begin{equation}
		X_2 - X_1= \frac{(X_{21}+X_{22})-(X_{11}+X_{12})}{2} = \frac{(X_{21}-X_{11})+(X_{22}-X_{12})}{2}.
\end{equation}
If we choose $X_2-X_1 \sim \mathcal{N}(0,\sigma^2)$, then we may choose $X_{21}-X_{11} \sim X_{22}-X_{12} \sim \mathcal{N}(0,2\sigma^2)$.
Hence, from this follows the scaling by $1/h'$ in Equation (\ref{eqn:gaussdiffs}).

%
%
Similarly for the Cauchy differences we may write
\begin{equation}
\begin{split}
X_{j,j'} - X_{j-1,j'} &\sim \mathrm{Cauchy}(\lambda h)\cr
X_{j,j'} - X_{j,j'-1} &\sim \mathrm{Cauchy}(\lambda h').
\end{split}
\end{equation}
The result is obtained through similar argumentation as in the Gaussian case.
%
%
If we make the grid denser along the regularization direction, this corresponds again to one-dimensional Cauchy random walk, i.e.\ we simply scale with discretization step $h$.

In the second coordinate direction we take again the same division as in Figure \ref{fig:dens}.
	 If  $X_2-X_1 \sim \mathrm{Cauchy}(\lambda)$, then $X_{21}-X_{11} \sim X_{22}-X_{12} \sim  \mathrm{Cauchy}(\lambda)$  (note independent linear combination of Cauchy distributed random variable).
%
%

%
Through similar arguments, we can actually see that the scaling of the general L\'evy $\alpha$-stable difference prior is 
\begin{equation}
\begin{split}
X_{j,j'} - X_{j-1,j'} &\sim S_\alpha\left(h'^{(\alpha-1)/\alpha)}h^{1/\alpha}, \beta,0\right) \cr
X_{j,j'} - X_{j,j'-1} &\sim S_\alpha\left(h^{(\alpha-1)/\alpha)}h'^{1/\alpha}, \beta,0\right).
\end{split}
\end{equation}
This, naturally, implies similar constructions in the higher dimensions.

\begin{figure}[htp]
\begin{center}
  \subfigure[Realization of 2D Cauchy prior]{\includegraphics[width=0.49\textwidth]{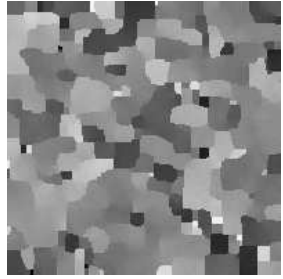}} \\
  \subfigure[Realization of 2D Gaussian prior]{\includegraphics[width=0.49\textwidth]{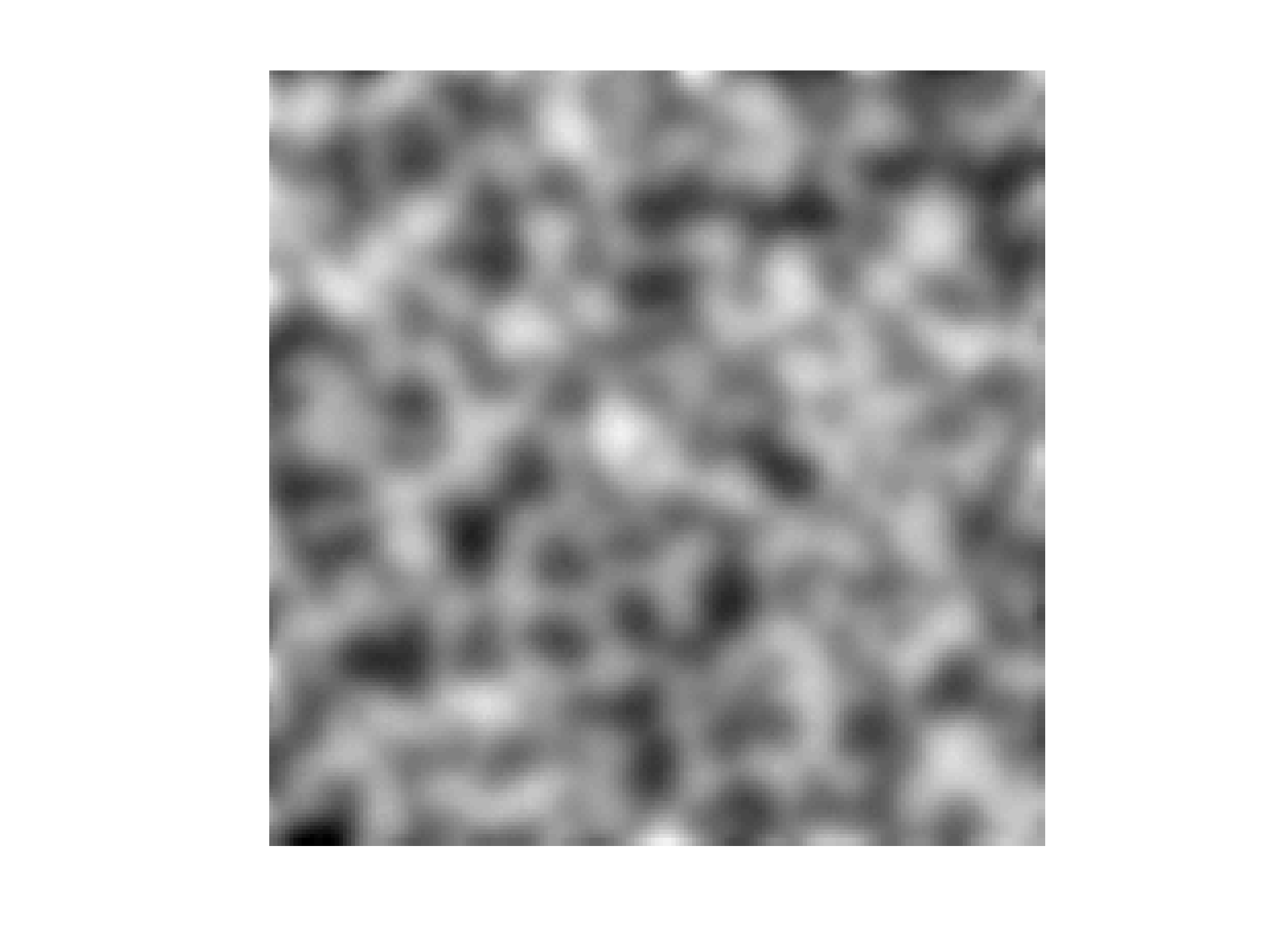}}
  \subfigure[Realization of 2D TV prior]{\includegraphics[width=0.49\textwidth]{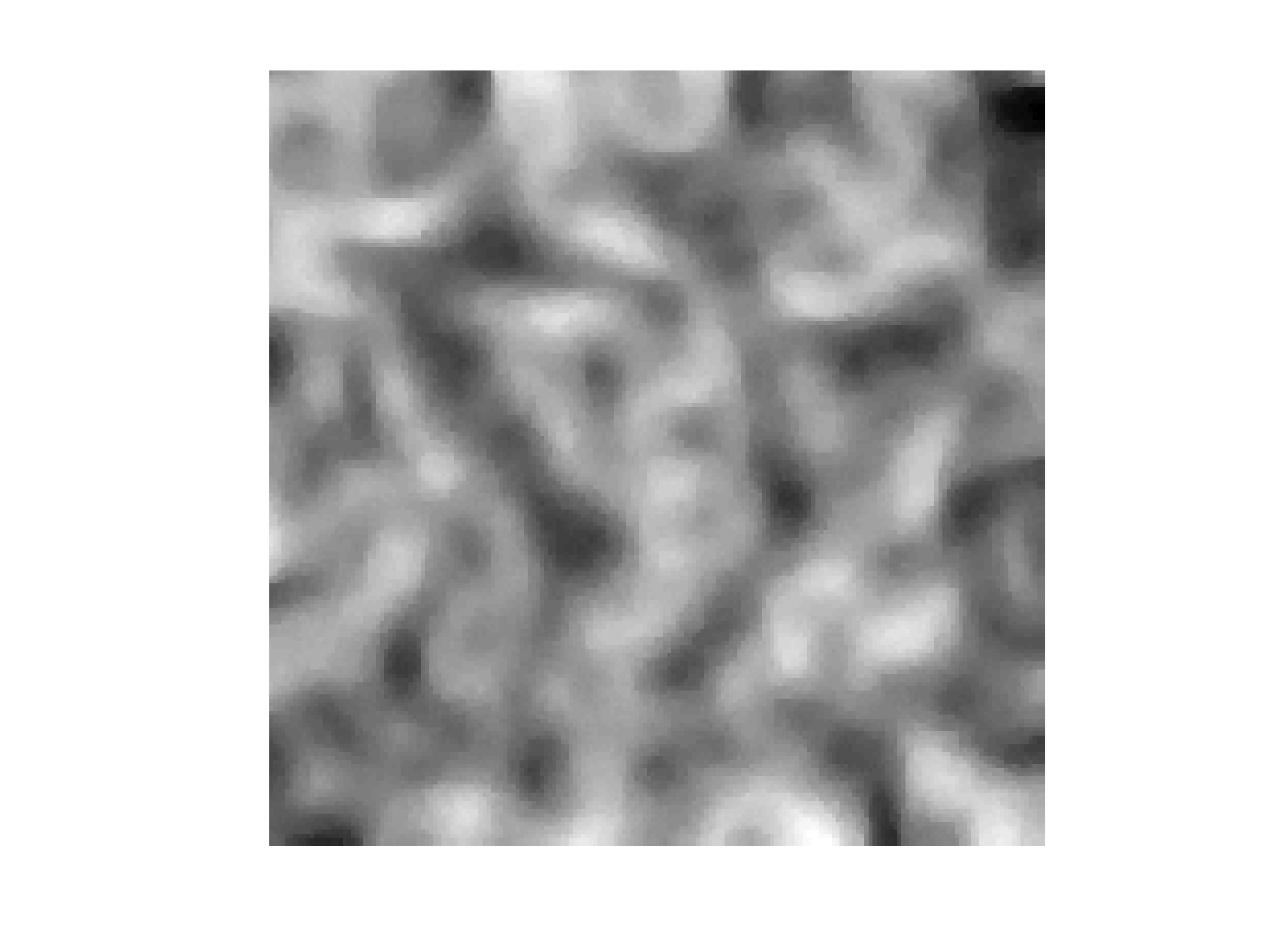}}

  

  \caption{Realizations of two-dimensional Cauchy, Gaussian and priors showing blocky and smooth structures.
  } \label{fig:2d_realisations}
  \end{center}
\end{figure}

In Figure \ref{fig:2d_realisations}, we have plotted realizations of two-dimensional Cauchy, Gaussian  and TV priors.
We have used zero-boundary conditions, but we have cropped the image in order to demonstrate blocky or smooth structures of the realizations. 
A notable feature in the Cauchy prior realization is that the realization does not look isotropic.
There seems to be more edges in vertical and horizontal directions than in other directions.
%
%
%
We may ask if this is a problem?
We suppose not, but this should be studied in more detail in subsequent papers. 
If we consider the edge-preserving features of the two-dimensional Cauchy prior, we can study processes along either of the coordinate directions.
These one-dimensional processes look practically similar to the  one-dimensional Cauchy random walks.
For comparison purposes, we have also plotted realization of a Gaussian prior obtained in a similar manner, which shows that the prior promotes smoothness, not jumps.
The Gaussian prior is also isotropic.
Realization of a TV prior is not as blocky as the Cauchy prior, but resembles more  the Gaussian prior.
However, the TV prior realization is not isotropic, as we use $\vert X_{j,j'}-X_{j-1,j'}\vert +\vert X_{j,j'}-X_{j,j'-1}\vert$, inside the exponential function of the TV prior, and we should use $\sqrt{\vert X_{j,j'}-X_{j-1,j'}\vert^2 +\vert X_{j,j'}-X_{j,j'-1}\vert^2}$  for  isotropic TV prior.


\subsection{Two-dimensional X-ray tomography}

\begin{figure}[htp]
\begin{center}
	\includegraphics[width=3in]{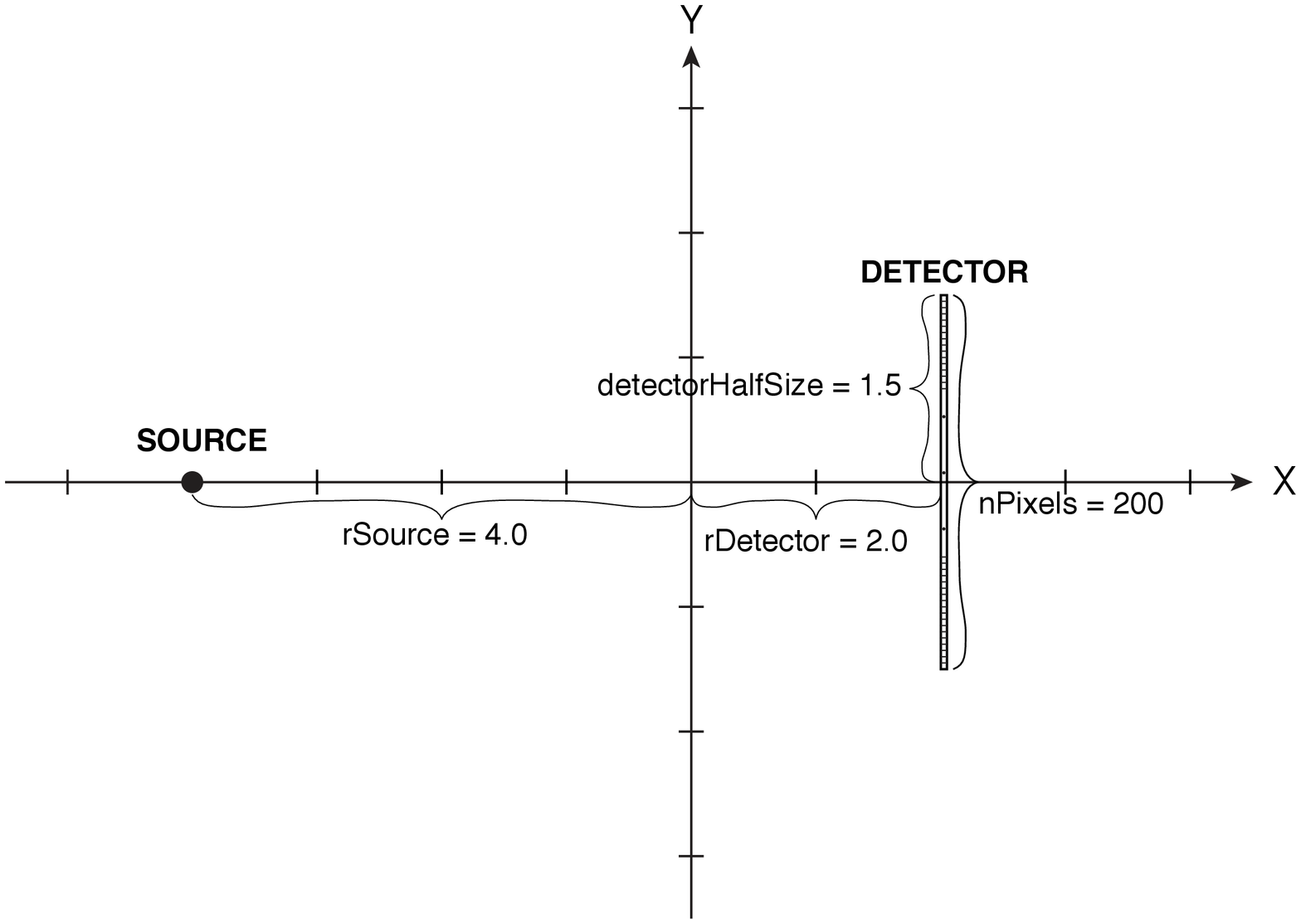}
	\caption{Measurement geometry for the synthetic two-dimensional tomography example.}
	\label{fig:tomogeom}
\end{center}
\end{figure}

Now we will apply the developed method to X-ray tomography.
We take the Shepp-Logan phantom to be the original unknown.
The measurements correspond to two-dimensional fan-beam geometry such that the distance between the source and the detector to the rotation axis is 4
and 2, respectively. 
The width of the detector is 3 and the detector element contains 200 pixels.  
In Figure \ref{fig:tomogeom}, we have depicted the geometry of the measurements. 

We make measurements at 41 different angles (from $-10^\circ$ to $190^\circ$ with $5^\circ$ steps).
We add white measurement noise with approximately 0.1\% noise level.
In order to get convergence of the MCMC chain, we use 5 million as the chain length, and 2.5 million as burn-in period. In order to save memory, we use only every 10th sample in the ensemble mean.
The computations are carried with a PC computational server with dual Xeon E5-2680 CPUs and 128 GB memory. 
The core parts of the algorithm are implemented using the MEX/C interface. 
Computation time for a reconstruction with $200 \times 200$ resolution is approximately 24 hours and, for a $400\times 400$ resolution, the time is approximately 96 hours. 
This indicates that the algorithm's computational cost is (nearly) proportional to the number of unknowns, as expected.


\begin{figure}[htp]
\begin{center}
  \subfigure[Unknown: the Shepp-Logan phantom]{\includegraphics[width=0.4\textwidth]{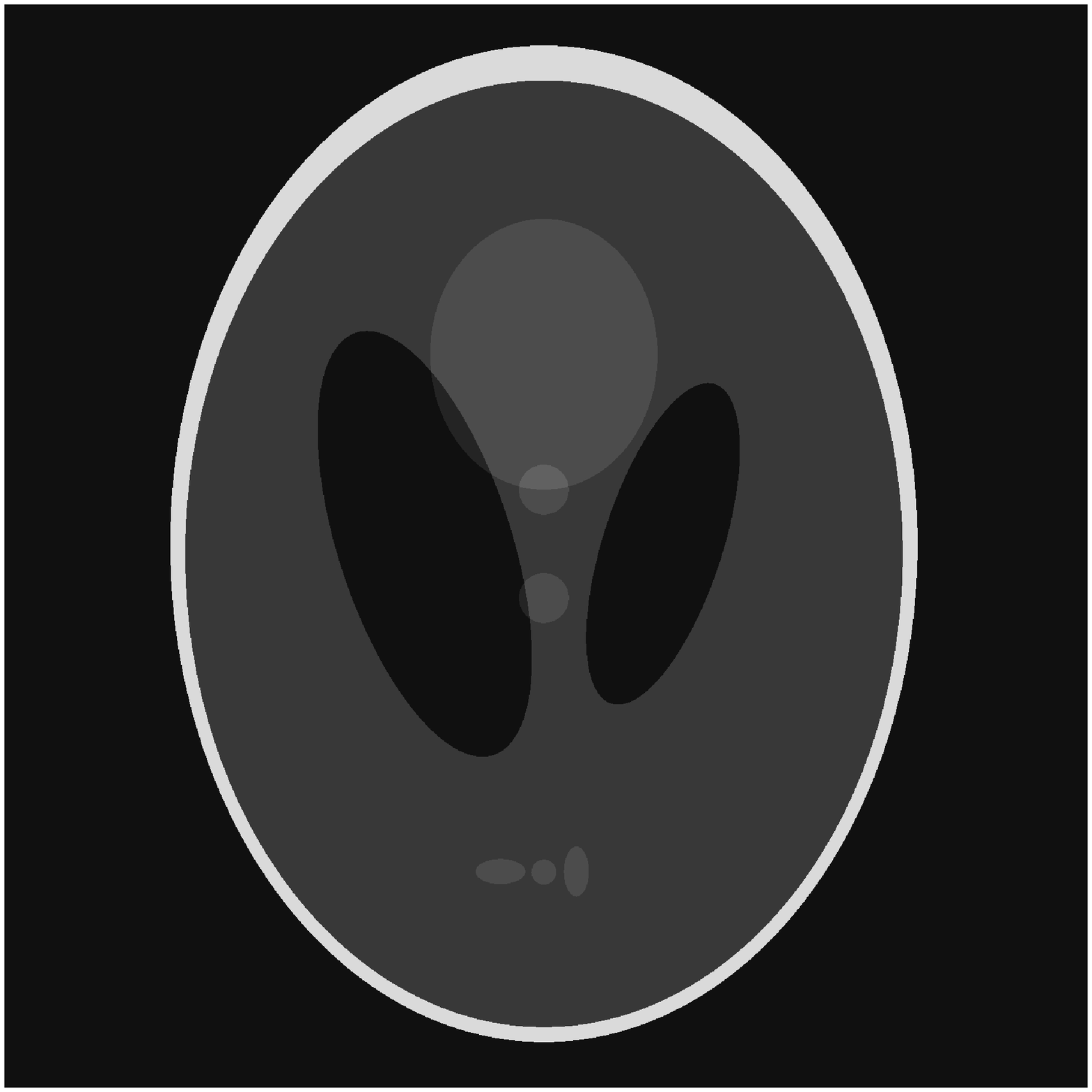}}
  \subfigure[Filtered back-projection estimate]{\includegraphics[width=0.4\textwidth]{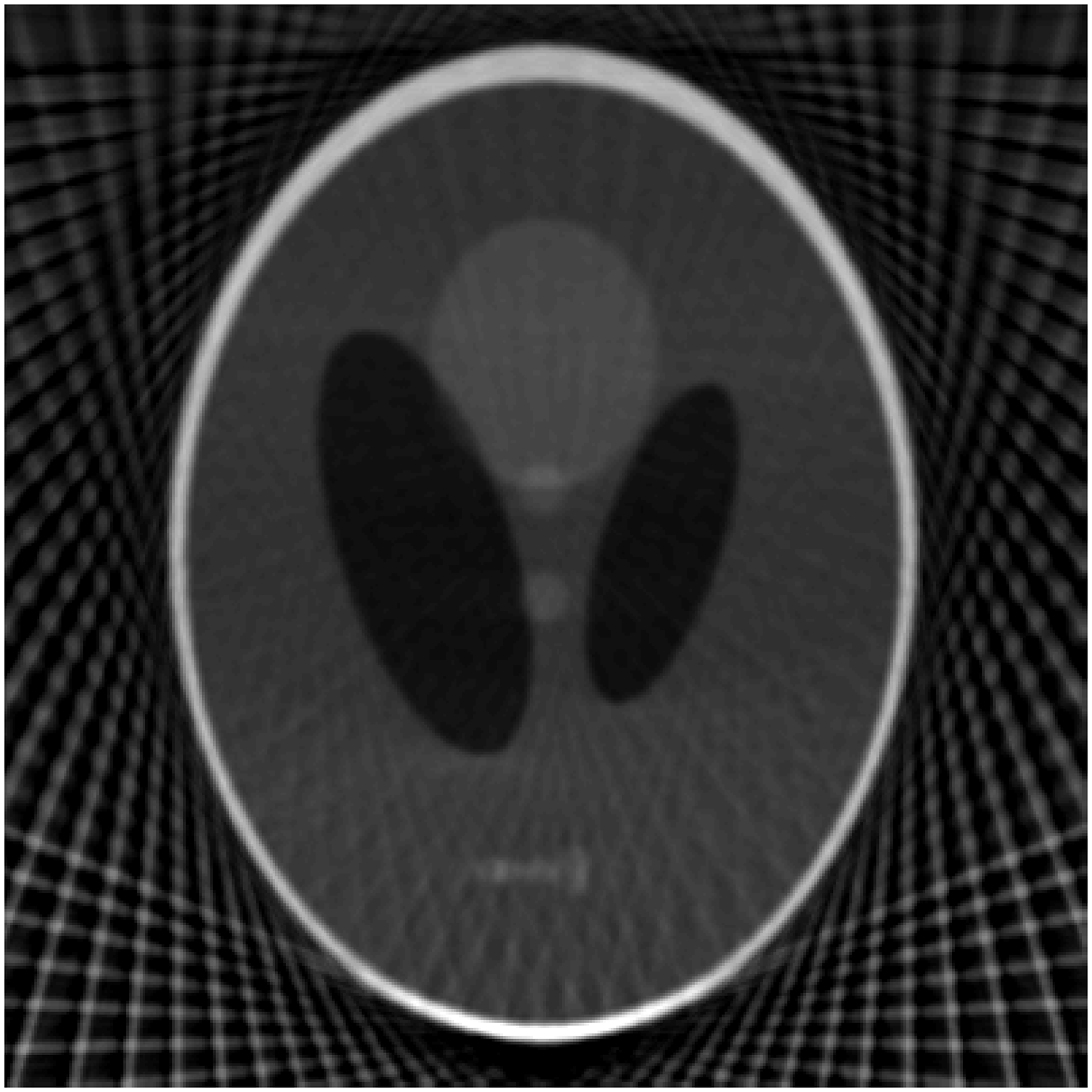}}

  \subfigure[MAP estimate with Cauchy prior]{\includegraphics[width=0.4\textwidth]{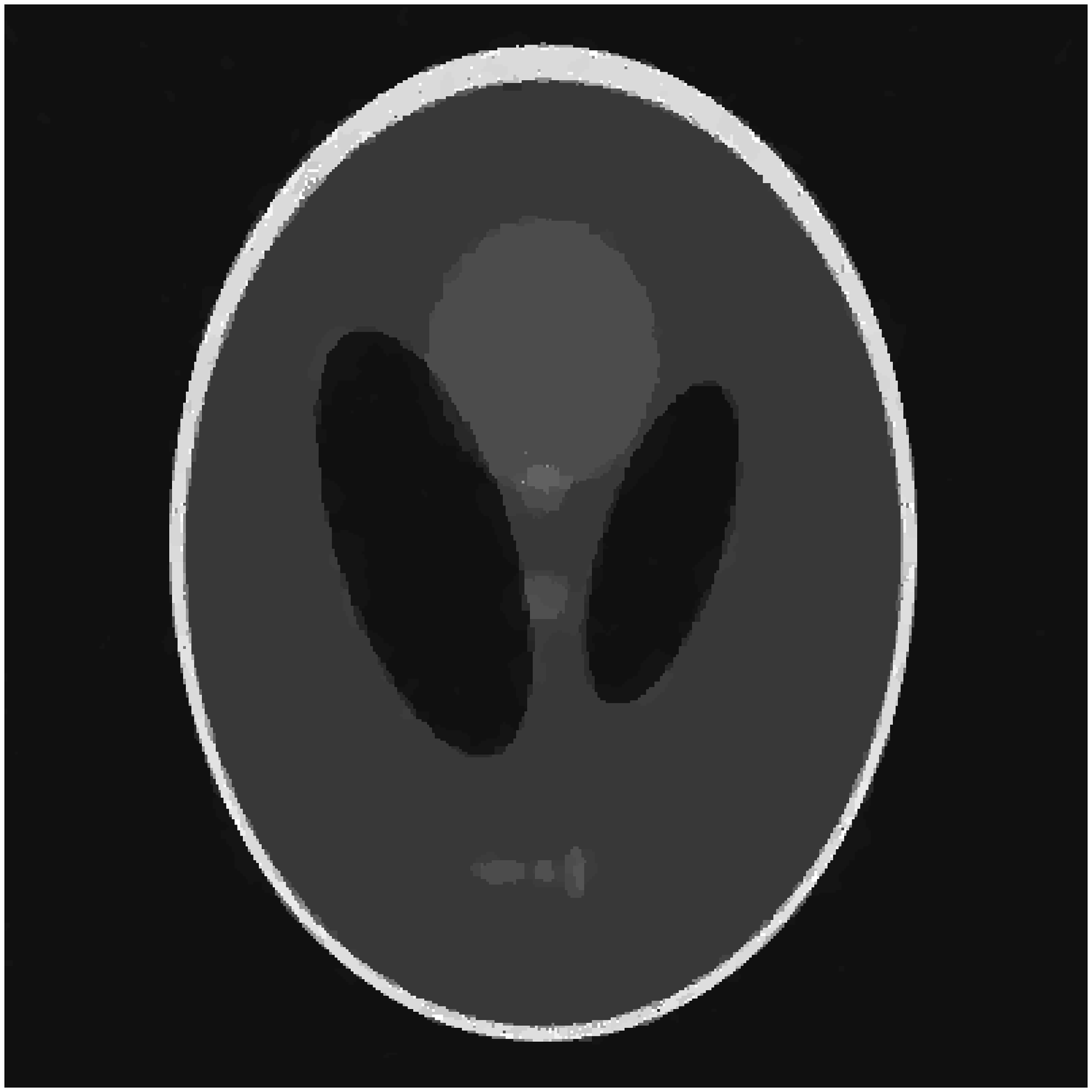}}
  \subfigure[CM estimate with Cauchy prior]{\includegraphics[width=0.4\textwidth]{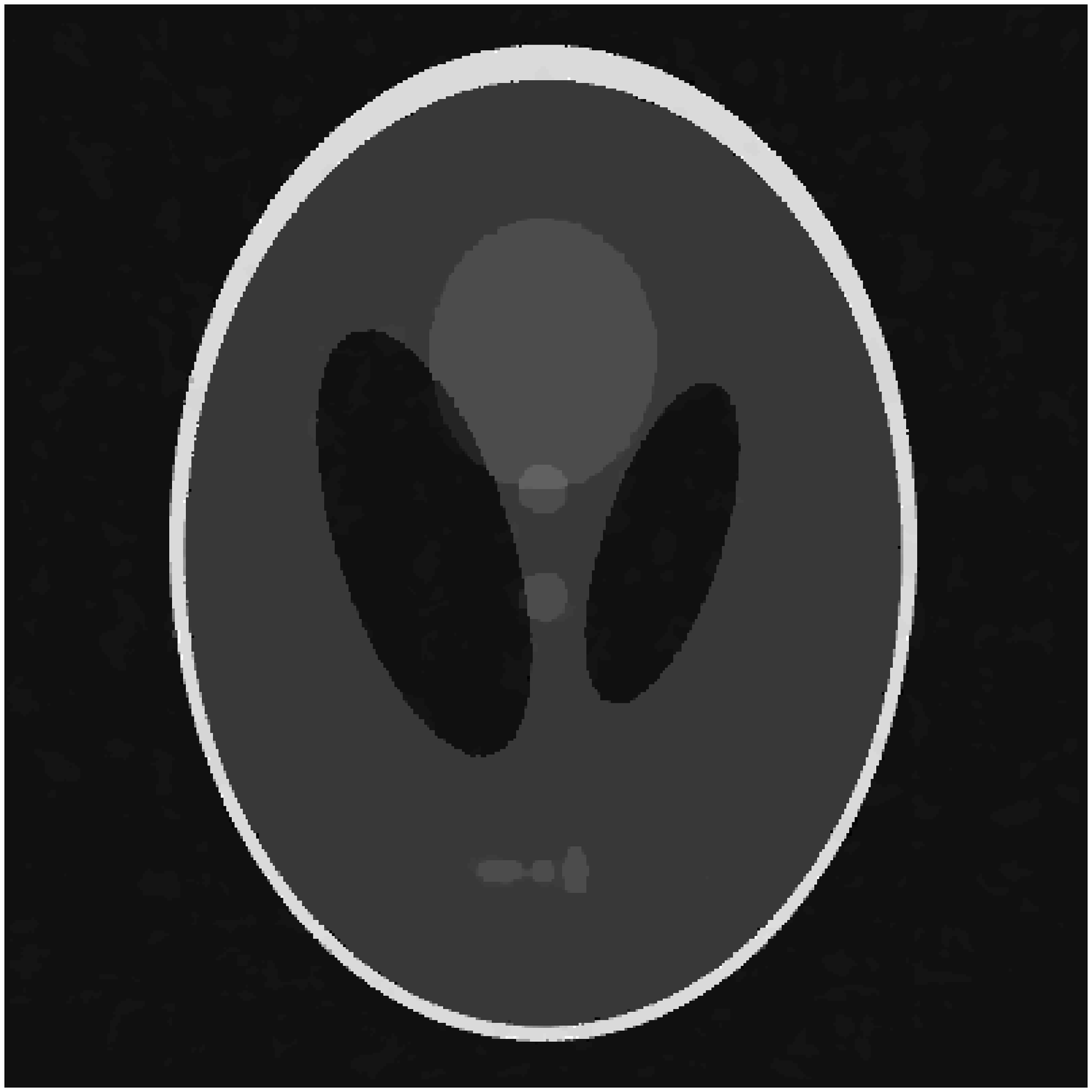}}

  \subfigure[CM estimate with TV prior]{\includegraphics[width=0.4\textwidth]{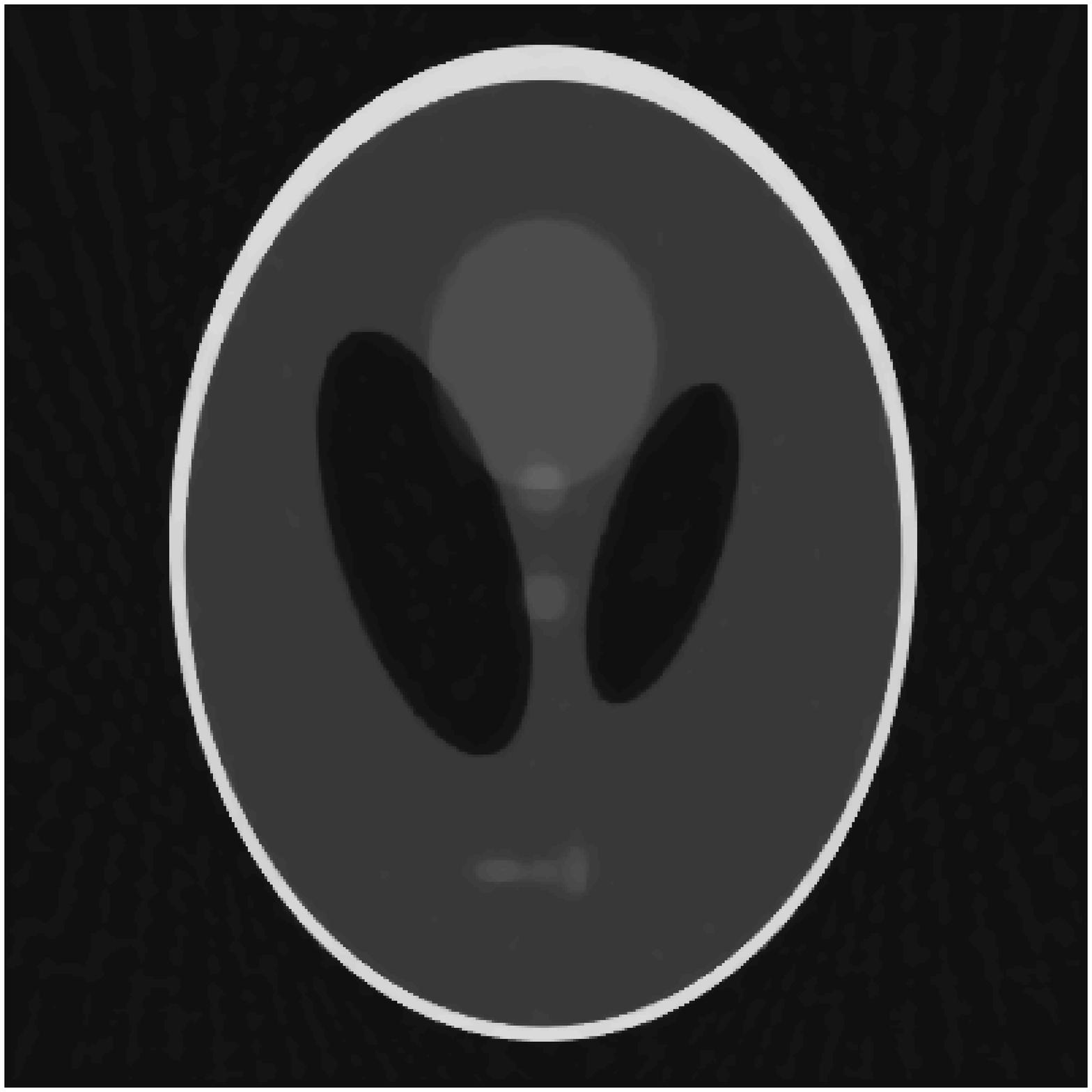}}
  \subfigure[CM estimate with Gaussian prior]{\includegraphics[width=0.4\textwidth]{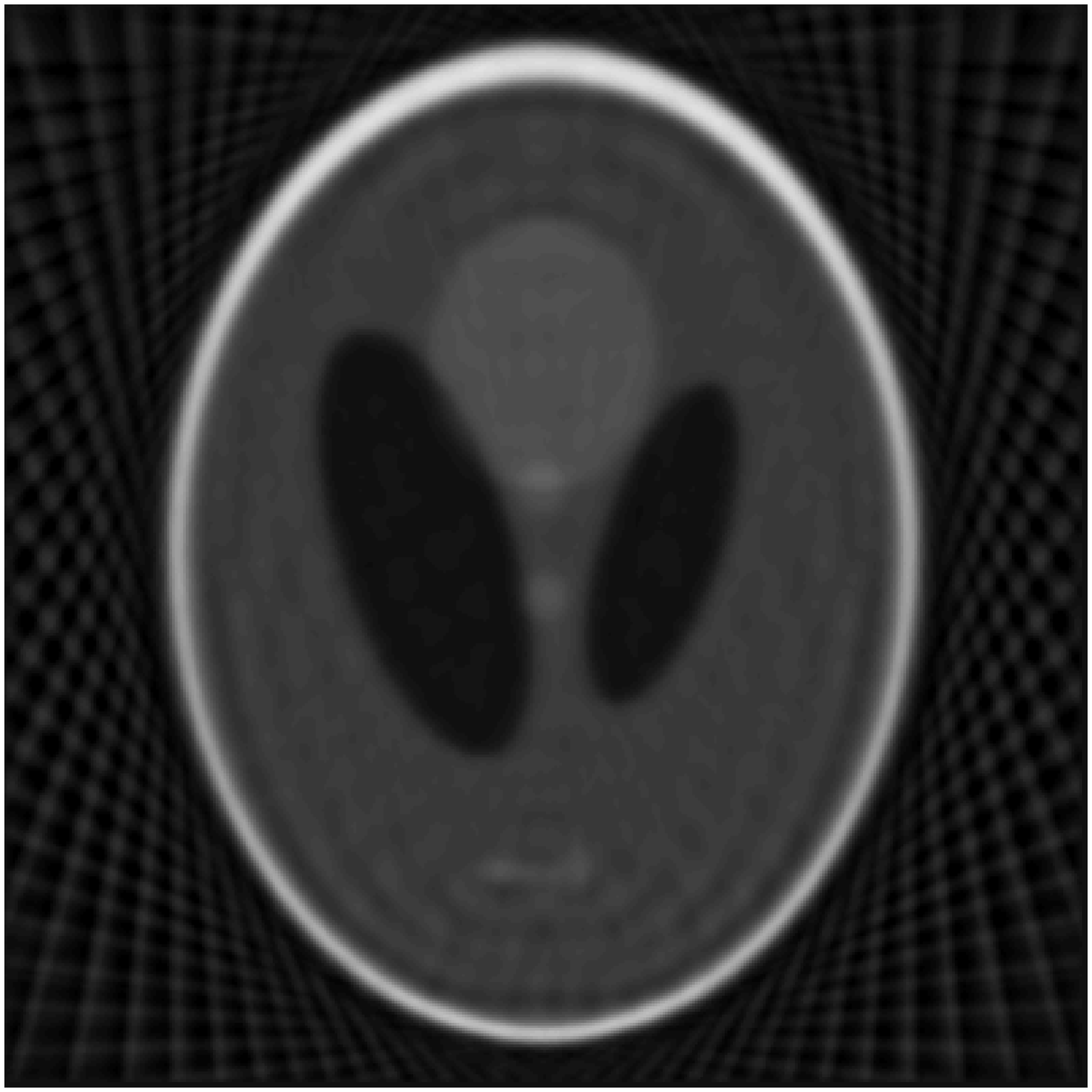}}

  \caption{Comparison of different estimates on lattice sized $400\times 400$ pixels.} \label{fig:comparison_tomo_priors}
  \end{center}
\end{figure}

In Figure \ref{fig:comparison_tomo_priors}, we have plotted the original unknown and reconstructions with five different methods.
The filtered back-projection estimate 
 recovers most of the features, but it has severe artifacts.
For the Cauchy prior, we have two different estimates, the MAP and CM estimates obtained with Gauss-Newton optimization methods and single-component Metropolis-Hastings, respectively.
Both methods capture the features well.
According to our experiments, the benefit of the MAP estimate is lower computational time, while the CM estimate has less artifacts in the reconstructions.
MAP estimates may have e.g.\ peaks in the reconstructions, which are not visible in the CM estimates.
The CM estimate with TV prior captures also the features well, but it is not as sharp as the Cauchy prior estimates. 
Also, the Cauchy estimate seems to be more robust, when we decrease the number of measurements or increase the noise.
The estimate with Gaussian prior recovers the main structures, but it has both bad artifacts and also smooth structures.



  


\begin{figure}[htp]
\begin{center}
  \subfigure[$200\times 200$]{\includegraphics[width=0.32\textwidth]{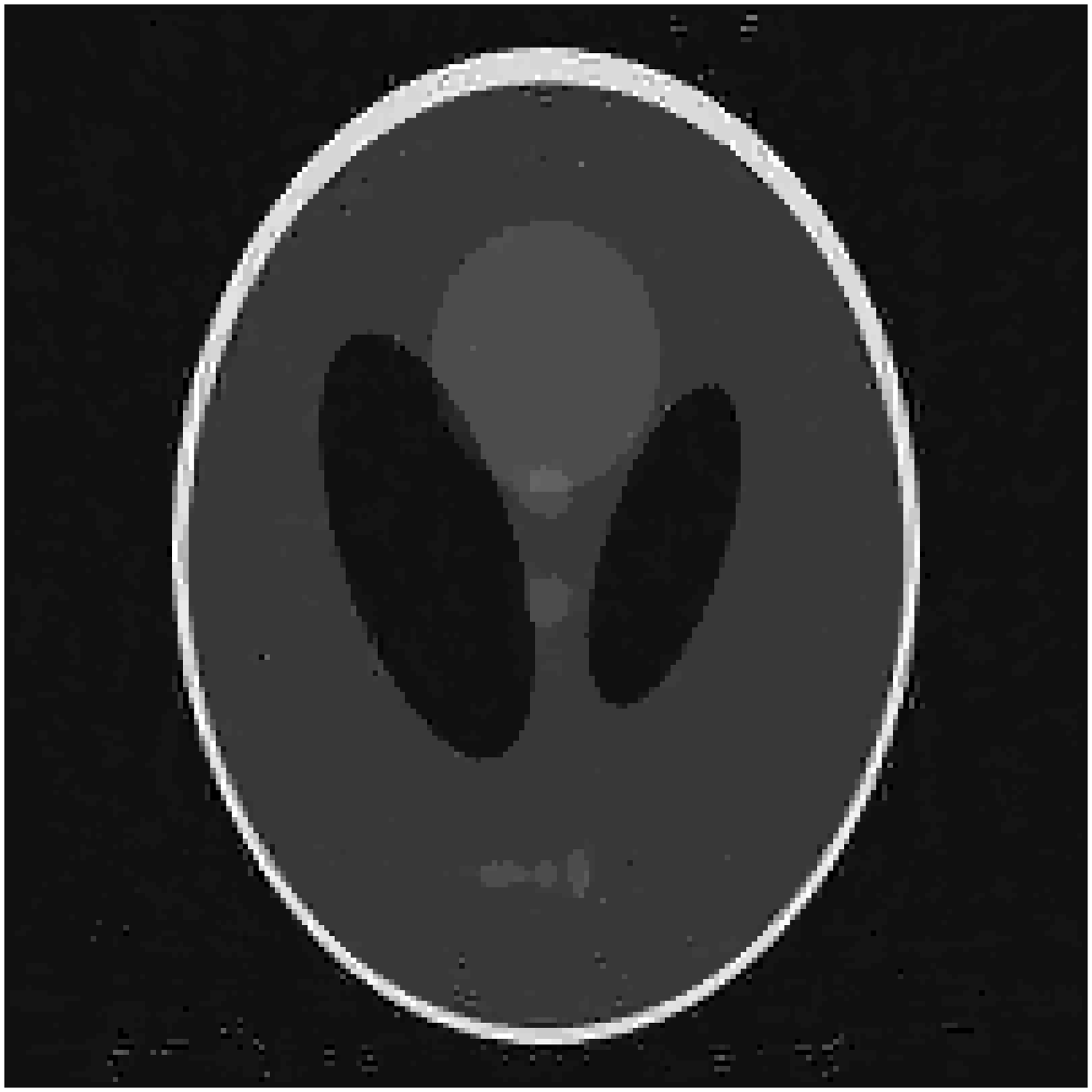}}
  \subfigure[$400\times 400$]{\includegraphics[width=0.32\textwidth]{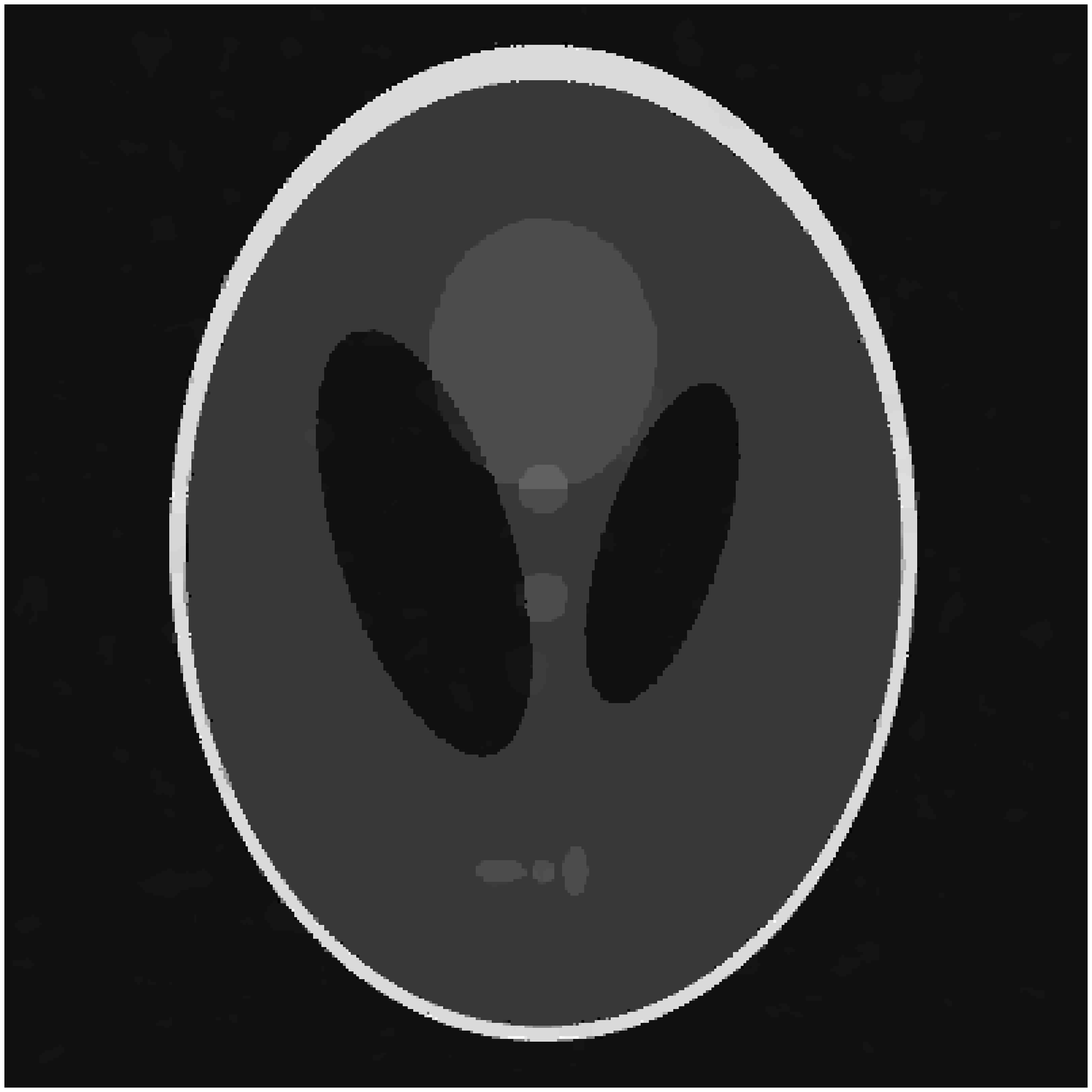}}
  \subfigure[$800\times 800$]{\includegraphics[width=0.32\textwidth]{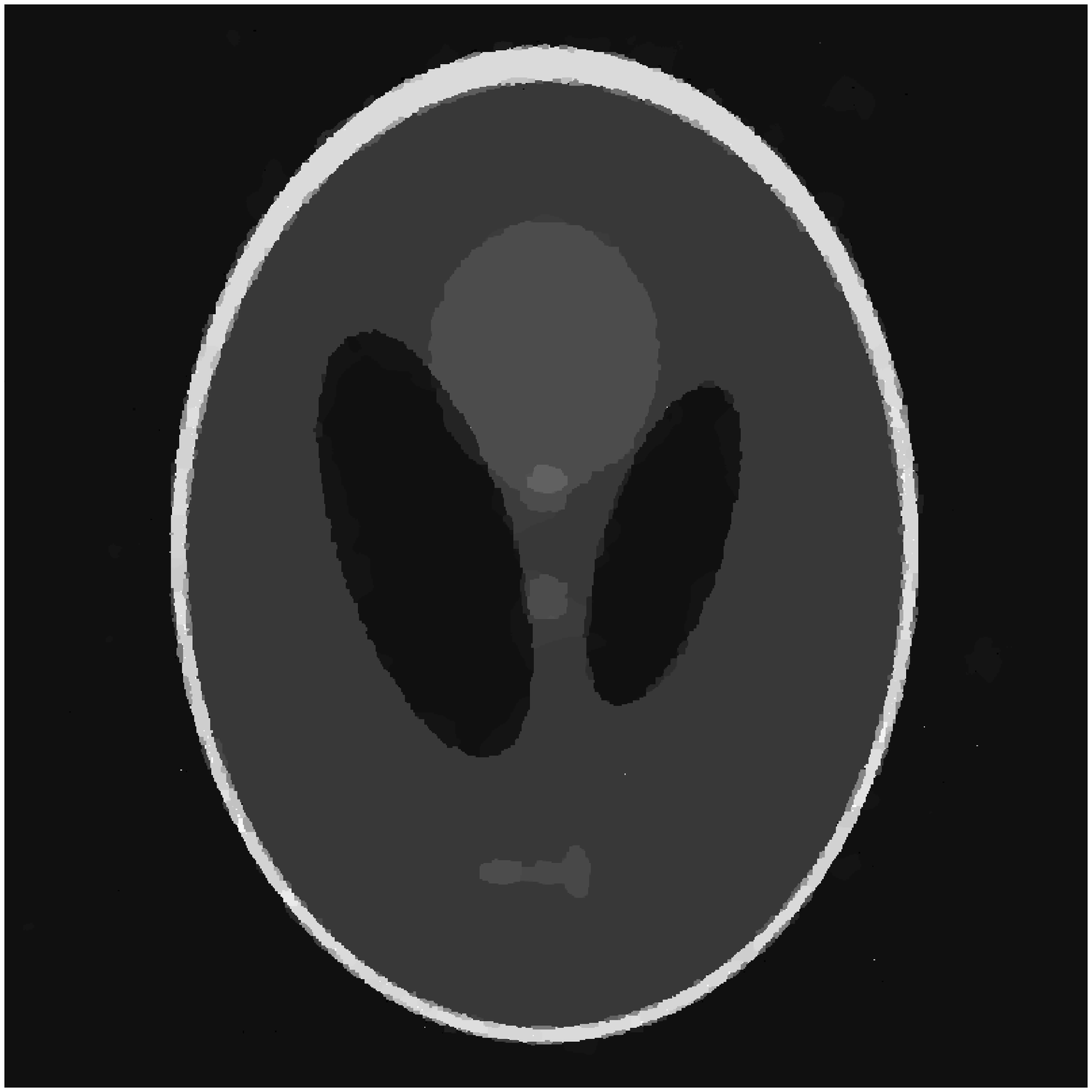}}

  \caption{Conditional mean estimates of the X-ray tomography problem on three different lattices.}¤
     \label{fig:comparison_tomo_priors_bb}
  \end{center}
\end{figure}

\begin{figure}[htp]
\begin{center}
  \subfigure[Horizontal]{\includegraphics[width=0.44\textwidth]{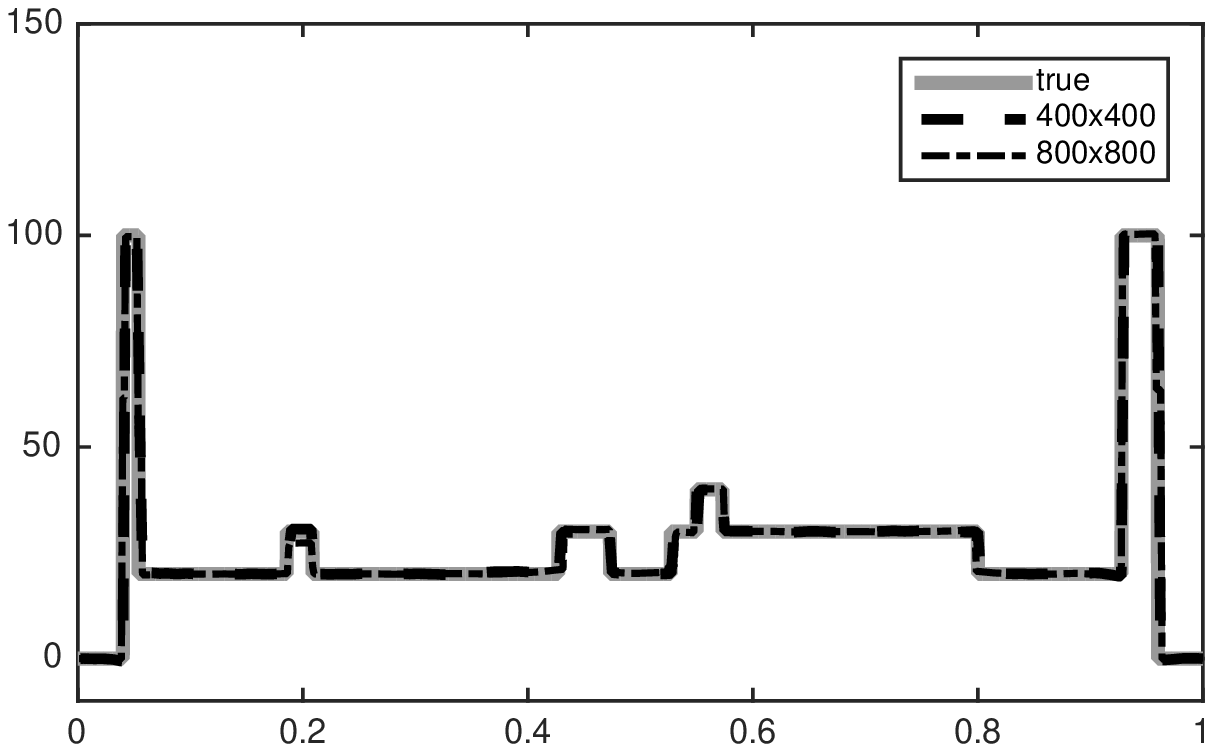}}
  \subfigure[Horizontal zoom-in]{\includegraphics[width=0.44\textwidth]{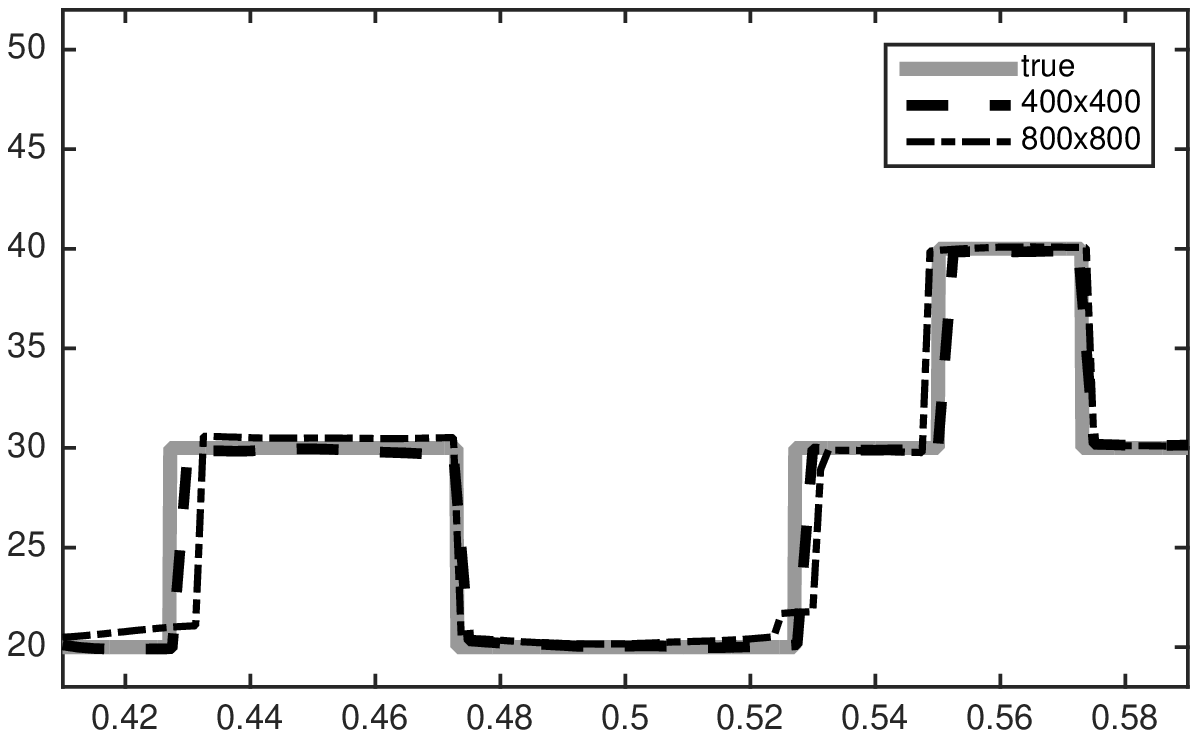}}

  \subfigure[Vertical]{\includegraphics[width=0.44\textwidth]{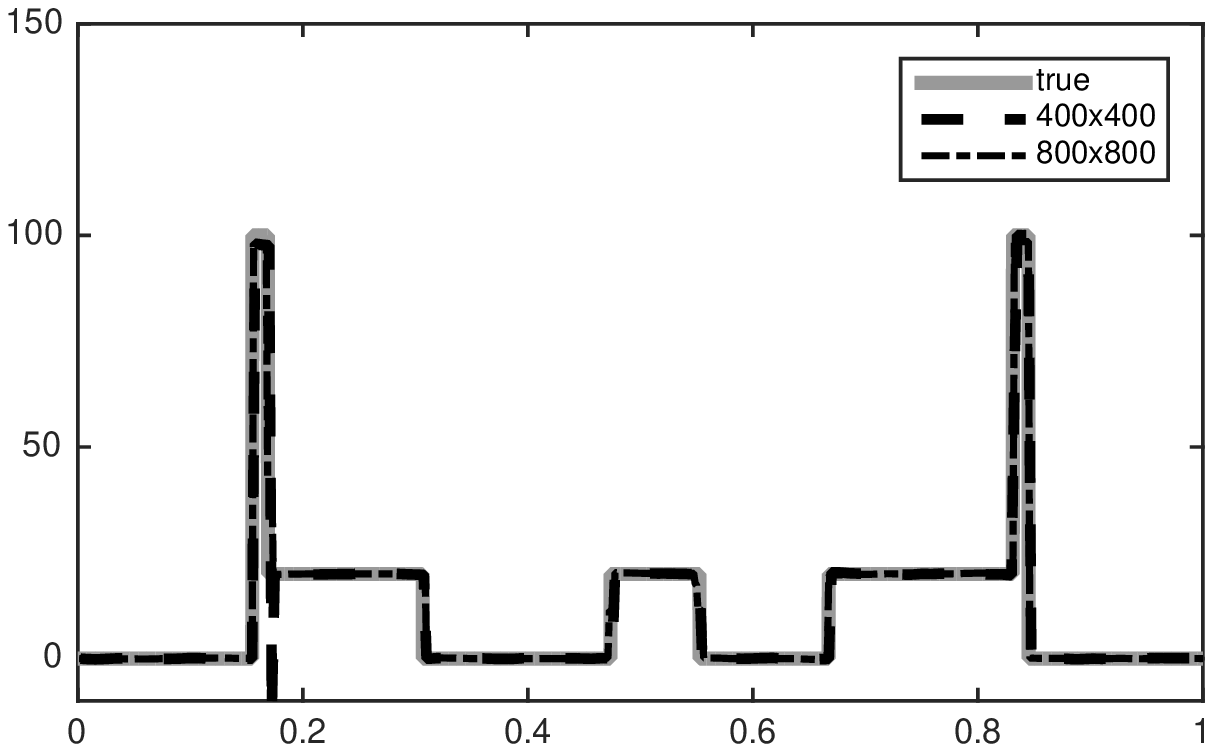}}
  \subfigure[Vertical zoom-in]{\includegraphics[width=0.44\textwidth]{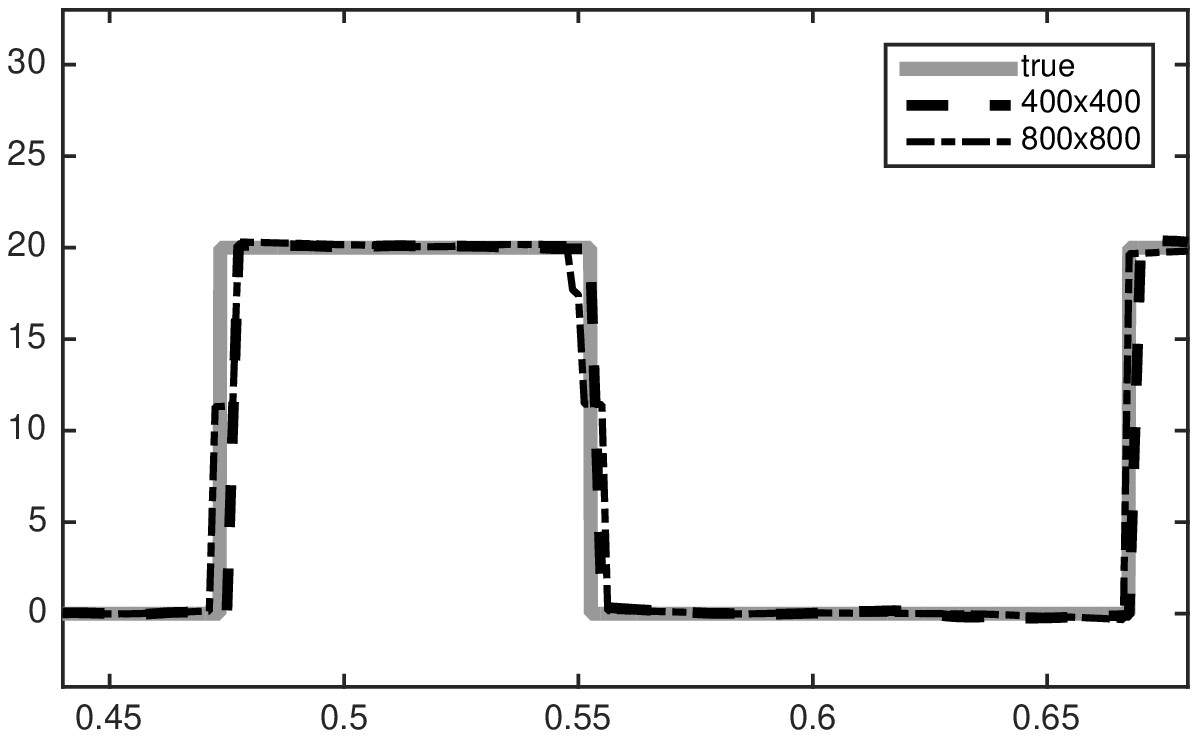}}
  \caption{Cross-sections of the conditional mean estimates from the middle of the computational domain.}
    \label{fig:comparison_tomo_priors_big}
  \end{center}
\end{figure}

In Figures \ref{fig:comparison_tomo_priors_bb} and \ref{fig:comparison_tomo_priors_big}, we have CM estimates on different size lattices.
The discretization has been taken into account as described earlier in this section, and hence the reconstructions look essentially same.

%

\section{Conclusion and discussion}
\label{sec:conclusion}

We have considered the construction of edge-preserving Bayesian inversion algorithms with Cauchy priors.
We drew estimators with  single-component Metropolis-Hastings.
When compared to Bayesian inversion with total variation prior, the methodology proposed promotes edge-preserving Bayesian inversion and gives good reconstructions for the tomography problem also when we have low number of measurements, the measurement noise starts to grow or when we have dense computational mesh.

In the  future studies, we should consider the discretization-invariance with mathematical rigor. 
Also, an exhaustive numerical study of L\'evy $\alpha$-stable random walks should be made.
%
While the Gaussian random walk with $\alpha=2$ is obviously some kind of a high end of stable distributions, the Cauchy random with $\alpha=1$ 
%
is in the middle of values $0 < \alpha \leq 2$. 
For $0< \alpha<1$ we probably have even more constant behavior between jumps, than we have  for the Cauchy prior, whereas  the cases for $1\leq \alpha\leq 2$ between Cauchy and Gaussian should provide a continuity from smoothing inference to edge-preserving inference.
%
%
Also, we could try more general prior models through stochastic partial differential of form $\mathcal{PX}=\mathcal{W}$, where $\mathcal P$ is some linear differential operator, $\mathcal X$ is unknown and $\mathcal W$ L\'evy noise.
This is similar to e.g.\ Mat\'ern field, i.e.\ a Gaussian Markov random field construction, except that we would have general L\'evy noise instead of the special case $\alpha=2$ Gaussian noise.

\section*{Acknowledgments}

This work has been funded by Engineering and Physical Sciences Research Council, United Kingdom (EPSRC Reference:	EP/K034154/1 --  Enabling Quantification of Uncertainty for Large-Scale Inverse Problems), European Research Council (ERC Advanced Grant 267700 - Inverse Problems) and Academy of Finland (application number 250215, Finnish Programme for Centre of Excellence in Research 2012-2017).

\end{document}